\documentclass[3pt]{elsarticle}
\pdfoutput=1
\usepackage{amsthm}
\usepackage[intlimits]{amsmath}
\usepackage{amssymb}
\usepackage{amsthm}
\usepackage{amsfonts}
\usepackage{bm}
\usepackage{mathdots}
\usepackage[normalem]{ulem}

\usepackage[utf8]{inputenc} 
\usepackage[T1]{fontenc}
\usepackage{times}
\usepackage{mathptmx}  
\DeclareSymbolFont{cmletters}{OML}{cmr}{m}{n}
\DeclareMathAlphabet{\mathcal}{OMS}{cmsy}{m}{n} 

\journal{arXiv}
\ifdefined\copyReadyOn
\usepackage{ecrc}
\volume{00}
\firstpage{1}
\runauth{G. Beer, Ch. Duenser, V. Mallardo}
\journalname{\journal}
\jid{IJNME}
\jnltitlelogo{IJNME}
\CopyrightLine{2015}{Published by Elsevier Ltd.}
\fi

\usepackage{tikz}
\usepackage{tikz-3dplot}
\usepackage{pgfplots}
\usepackage{pgfplotstable}
\usepackage{tikzscale}
\usepackage{standalone} 
\usepackage{hf-tikz}

\usepackage{forloop}
\usepackage{arrayjobx}

\usepackage{algorithmic}
\usepackage{algorithm}

\usepackage{booktabs} 
\usepackage{calc}
\usepackage{xcolor,pdfcomment,soul} 
\usepackage{xspace} 
\usepackage{multirow}

\usepackage{scrlayer-scrpage}
\usepackage{setspace}

\usepackage[labelformat=simple]{subcaption}   	
\usepackage{pdfpages} 

\usepackage{tabto} 

\usepackage{makeidx}
\makeindex

\usepackage{colortbl}

\usepackage{enumerate}

\usepackage{overpic}

\usepackage{listings} 

\usepackage[intlimits]{amsmath}
\usepackage{amssymb} 
\usepackage{amsthm}
\usepackage{amsfonts,mathrsfs}
\usepackage{bm}
\usepackage{mathtools} 

\usepackage{siunitx} 

\usepackage{xfrac}
\usepackage{listings}
\usepackage{color}
\definecolor{dkgreen}{rgb}{0,0.6,0}
\definecolor{gray}{rgb}{0.5,0.5,0.5}
\definecolor{mauve}{rgb}{0.58,0,0.82}

\usepackage[utf8]{inputenc} 
\usepackage[german,english]{babel}  
\usepackage[T1]{fontenc} 
\usepackage{lmodern}

\usepackage{geometry}
\usepackage{overpic}
\graphicspath{{pics/}{tmp/}}

\newboolean{inGray}
\setboolean{inGray}{false}

\newcommand{\addtoindex}[2][]{
    \ifthenelse { \equal{#1}{} }
    {#2\index{#2}\xspace}%
    {#2\index{#1}\xspace}%
}

\newcommand{\myVec}[1]{\mathbf{#1}}
\newcommand{\myVecGreek}[1]{\boldsymbol{#1}}

\newcommand\tens[2]{\mathsf{#1}} %

\newcommand{\NURBS}{R} 			
\newcommand\uu{\xi} 			
\newcommand\vv{\eta} 			
\newcommand{\domain}{\Omega}

\newcommand{\boundary}{\Gamma}

\newcommand\primary{u} 			
\newcommand\dual{t} 			

\providecommand\url[1]{\emph{#1}}

\newcommand\fund[1]{\tens{#1}{2}}

\newcommand\pt[1]{\boldsymbol{#1}}
\newcommand\sourcept{\tilde{\pt{x}}}
\newcommand\fieldpt{\hat{\pt{x}}}

\newenvironment{mytable}[4]%
{
  \begin{table}[#1]
    \def\mycap{#2}
    \def\mylabel{#3}
    \centering
    \begin{tabular}{#4}
      \toprule
}
{
  \bottomrule
  \end{tabular}
  \caption{\mycap}
  \label{\mylabel}
  \end{table}
}
\newcommand{\mytableheader}[1]{
       #1 \\ \midrule
}

{
  \begin{algorithm}
    \caption{#1}
    \label{#2}
    \begin{algorithmic}[1]
}   
{
\end{algorithmic}
\end{algorithm}
}

\newtheoremstyle{myremark}
{3pt}
{3pt}
{}
{}
{\itshape}
{:}
{.5em}
{}
\theoremstyle{myremark}
\newtheorem*{remark}{Remark}

\newcommand{\myeqref}[1]{equation~(\ref{#1})}	
\newcommand{\myfigref}[1]{Figure~\ref{#1}}

\newcounter{footnoteNumber} 

\DeclareCaptionFormat{listing}{#1#2#3}
\newcommand{\myalignatsinglelabel}[3]%
{%
    \begin{equation}
        \label{#2}
        \begin{alignedat}{#1} 
            #3
        \end{alignedat}
    \end{equation}
}

\newcommand{\myalignat}[2]%
{%
  \begin{alignat}{#1} 
    #2
   \end{alignat}
}

\usepackage{pgf,pgffor,pgfmath}
\pgfplotsset{compat=1.10}

\usetikzlibrary{external}
\usetikzlibrary{positioning,calc,intersections,arrows,3d}
\usetikzlibrary{shapes,snakes,fit,spy}
\usepgfplotslibrary{patchplots}
\usetikzlibrary{matrix}

\usepgfplotslibrary{groupplots}
\usepgfplotslibrary{fillbetween}

\def\pgfplotfontsizetitle{\small}
\def\pgfplotfontsizelegend{\small}
\def\pgfplotfontsize{\small}
\def\pgfplotfontsizetiny{\scriptsize}

\def\tikzfontsizetiny{\scriptsize}

\pgfplotsset{
  mystyle/.style ={%
    grid = major,
    every tick label/.append style={font=\pgfplotfontsizetiny},
    every axis label/.append style={font=\pgfplotfontsize},
    legend style={font=\pgfplotfontsizelegend},
    label style={font=\pgfplotfontsize},
    title style={font=\pgfplotfontsizetitle},
    /pgf/number format/set thousands separator = {}, 
  }
}%

\makeatletter
\pgfplotsset{
    myIgnoreRowModulo2/.style args={#1}{%
        /pgfplots/x filter/.code={%
        \let\xValue\pgfmathresult 
        \pgfmathparse{int(mod(int(\coordindex),int(2))} \pgfmathresult 
        \ifnum#1=\pgfmathresult
            \def\pgfmathresult{} 
        \else
            \pgfmathparse{\xValue} \pgfmathresult 
        \fi
        }
    } 
}
\makeatother

\colorlet{drawblue}      {blue!80!white}
\colorlet{drawred}       {red!80!white}
\colorlet{drawgray}      {gray}
\definecolor{drawgreen}  {RGB}{44,162,95}
\colorlet{drawpurple}    {purple}
\colorlet{draworange}    {orange}
\colorlet{drawlime}      {lime!80!black}
\colorlet{drawartichoke} {yellow!60!black}
\colorlet{TUGgray}{black!15}
\definecolor{TUGred}{RGB}{247,1,70}
\definecolor{IFBblue}{RGB}{51,112,169}

\definecolor{basisColor1}{RGB}{59,76,192}
\definecolor{basisColor2}{RGB}{87,117,225}
\definecolor{basisColor3}{RGB}{119,154,247}
\definecolor{basisColor4}{RGB}{152,185,255}
\definecolor{basisColor5}{RGB}{184,208,249}
\definecolor{basisColor6}{RGB}{195,209,230}	
\definecolor{basisColor7}{RGB}{247,200,190}	
\definecolor{basisColor8}{RGB}{247,187,160}
\definecolor{basisColor9}{RGB}{244,154,123}
\definecolor{basisColor10}{RGB}{229,112,88}
\definecolor{basisColor11}{RGB}{203,62,56}
\definecolor{basisColor12}{RGB}{180,4,38}

\definecolor{basisColor8sw}{RGB} {189,189,189}
\definecolor{basisColor9sw}{RGB} {150,150,150}
\definecolor{basisColor10sw}{RGB}{115,115,115}
\definecolor{basisColor11sw}{RGB}{91,91,91}
\definecolor{basisColor12sw}{RGB}{37,37,37}

\colorlet{myblue}    {blue}
\colorlet{myred}     {red}
\colorlet{mygreen}   {drawgreen}
\colorlet{mypurple}  {purple}
\colorlet{myorange}  {orange}

\ifthenelse{\boolean{inGray}}{ 

\pgfplotscreateplotcyclelist{mycyclelist}{
  basisColor12sw,semithick, every mark/.append style={solid, fill=.},mark=*           \\%
  basisColor11sw,semithick, every mark/.append style={solid, fill=.},mark=square*     \\%
  basisColor10sw,semithick, every mark/.append style={solid, fill=.},mark=triangle*   \\%
  basisColor9sw , semithick,every mark/.append style={solid, fill=.},mark=diamond*    \\%
  basisColor8sw , semithick,every mark/.append style={solid, fill=.},mark=otimes    \\%
}
}{

\pgfplotscreateplotcyclelist{mycyclelist}{
  basisColor12,semithick,   every mark/.append style={solid, fill=.},mark=*           \\%
  basisColor11,semithick,   every mark/.append style={solid, fill=.},mark=square*     \\%
  basisColor10,semithick,   every mark/.append style={solid, fill=.},mark=triangle*   \\%
  basisColor9, semithick,   every mark/.append style={solid, fill=.},mark=diamond*    \\%
  basisColor8, semithick,   every mark/.append style={solid, fill=.},mark=otimes    \\%
}
}

\ifthenelse{\boolean{inGray}}{ 
\pgfplotscreateplotcyclelist{mycyclelistcompare}{
  basisColor12sw, semithick, loosely dashdotdotted     \\%
  basisColor11sw, semithick, loosely dashed      		\\%
  basisColor10sw, semithick, dashed \\
  basisColor9sw,  semithick, solid               	\\%
  basisColor12sw, only marks,   every mark/.append style={solid, fill=.},mark=otimes    \\%
  basisColor11sw, only marks,   every mark/.append style={solid, fill=.},mark=triangle*   \\%
  basisColor10sw, only marks,   every mark/.append style={solid, fill=.},mark=square*     \\%
  basisColor9sw,  only marks,   every mark/.append style={solid, fill=.},mark=*          \\%
}

}{
\pgfplotscreateplotcyclelist{mycyclelistcompare}{
  basisColor12, semithick, loosely dashdotdotted     \\%
  basisColor11, semithick, loosely dashed      		\\%
  basisColor10, semithick, dashed \\
  basisColor9,  semithick, solid               	\\%
  basisColor12, only marks,   every mark/.append style={solid, fill=.},mark=otimes    \\%
  basisColor11, only marks,   every mark/.append style={solid, fill=.},mark=triangle*   \\%
  basisColor10, only marks,   every mark/.append style={solid, fill=.},mark=square*     \\%
  basisColor9,  only marks,   every mark/.append style={solid, fill=.},mark=*          \\%
}
}

\ifthenelse{\boolean{inGray}}{ 
\tikzset{mycyclelistcompareReferenceA/.style={basisColor12sw,solid}}
\tikzset{mycyclelistcompareTestA/.style={basisColor12sw,only marks,mark=otimes}}
}{
\tikzset{mycyclelistcompareReferenceA/.style={basisColor12,solid}}
\tikzset{mycyclelistcompareTestA/.style={basisColor12,only marks,mark=otimes}}
}

\tikzset{helpline/.style={thin,dashed}}
\tikzset{labelline/.style={thin}}
\tikzset{referencePath/.style={dotted,very thick}}
\tikzset{helparrow/.style={thin,arrows={-latex}}}
\tikzset{axis/.style={thin,arrows={->}}}
\tikzset{force/.style={thick,arrows={->}}}
\tikzset{forceInverse/.style={thick,arrows={<-}}}
\tikzset{Gamma/.style={ultra thick}}
\tikzset{controlPoly/.style={draw=black}}
\tikzset{GammaFill/.style={fill=lightgray,fill opacity=0.5}}
\tikzset{colorDiri/.style={drawgreen}}
\tikzset{GammaFillDiri/.style={fill=drawgreen,fill opacity=0.5}}
\tikzset{initialgrid/.style={thin,gray}}
\tikzset{addgridline/.style={dashed,gray}}
\tikzset{trimmingcurve/.style={thick}}
\tikzset{boundingbox/.style={thick, dotted}}
\tikzset{parameterSpace/.style={ }}
\tikzset{basisfunction/.style={very thick,smooth}}
\tikzset{bspline/.style={very thick,smooth,red}}
\tikzset{intersectioncurve/.style={dashed,thick}}
\tikzset{integrationRegionEdge/.style={dashed}}
\tikzset{pointer/.style={arrows={-latex}}}

\tikzstyle{anode}= [circle, inner sep=1.3pt, draw, fill=black]
\tikzstyle{gausspoint}=[shape=circle,draw=black,fill=black,inner sep=1.1pt]
\tikzstyle{controlPoint}=[shape=circle,draw=black,fill=black,thin,inner sep=0pt,minimum size=1.5mm]
\tikzstyle{abscissaPoint}=[shape=circle,draw=black,fill=white,thin,inner sep=0pt,minimum size=1.5mm]
\tikzstyle{anchorPoint}=[shape=circle,draw=black,fill=black,thin,inner sep=0pt,minimum size=1.5mm]
\tikzstyle{anchorPointDeg}=[shape=circle,draw=black,fill=TUGred,thin,inner sep=0pt,minimum size=1.5mm]
\tikzstyle{anchorPointDegD}=[shape=cross out,thick,draw=black,inner sep=0pt,minimum size=1.5mm]
\tikzstyle{trimmingIntersectionPoint}=[shape=cross out,thick,draw=black,inner sep=0pt,minimum size=1.5mm]

\tikzset{%
  highlight/.style={rectangle,rounded corners,fill=red!60,draw,fill opacity=0.125,thick,inner sep=0pt}
}

\usepackage{colortbl}

\def\trianglecolor{black}
\newcommand{\upperSlopeTriangle}[4] 	
{
	\def\trianglecolor{black}
	\addplot[forget plot, domain=#3:#4,color=\trianglecolor,samples=2]{  #2 / (x^#1) } node (A1) [pos=1] {}; 
	\addplot[forget plot, domain=#3:#4,color=\trianglecolor,samples=2]{   #2  / (#3^#1)} node (A2) [pos=1] {} node [anchor=south,pos=0.5,black] {\tikzfontsizetiny $1$};
	\draw[color=\trianglecolor] (A1.center) -- (A2.center) node [anchor=west,pos=0.5,black] {\tikzfontsizetiny #1};
}

\newcommand{\lowerSlopeTriangle}[4] 	
{
	\def\trianglecolor{black}
	\addplot[forget plot, domain=#3:#4,color=\trianglecolor,samples=2]{  #2 / (x^#1) } node (A1) [pos=0] {}; 
	\addplot[forget plot, domain=#3:#4,color=\trianglecolor,samples=2]{   #2  / (#4^#1)} node (A2) [pos=0] {} node [anchor=north,pos=0.5,black] {\tikzfontsizetiny $1$};
	\draw[color=\trianglecolor] (A1.center) -- (A2.center) node [anchor=east,pos=0.5,black] {\tikzfontsizetiny #1};
}

\newcommand{\myaddgraphic}[5]
{
 \node[anchor=south west,inner sep=0] (image) {\phantom{\includegraphics[#2]{#1}}};
  \begin{scope}[x={(image.south east)},y={(image.north west)}]
      
      \begin{scope}
          
          #5
          
          \node[anchor=south west,inner sep=0] {\includegraphics[#2]{#1}};
      \end{scope} 
      
      #4
      
      \pgfmathparse{int(#3)} \let\gridIndicator\pgfmathresult
      \ifthenelse{ \gridIndicator = 1 }
      {
          \draw[help lines,xstep=.1,ystep=.1] (0,0) grid (1.001,1.001);
          \foreach \x in {1,...,9} { \node [anchor=north] at (\x/10,0) {\x};}
          \foreach \y in {1,...,9} { \node [anchor=east] at (0,\y/10) {\y};}
      }{}
      
  \end{scope}    
}

\tikzstyle{reverseclip}=[insert path={(current page.north east) --
    (current page.south east) --
    (current page.south west) --
    (current page.north west) --
    (current page.north east)}
]

\newcounter{itR}
\newcommand{\bsplinevalue}[5] 
{                    				
	\newarray\vKnots
	\newarray\vN
	\newarray\vNumeratorL
	\newarray\vNumeratorR
	\newarray\vSave

	\readarray{vKnots}{#1}
	\readarray{vN}{1}
	\readarray{vSave}{0}
	\readarray{vNumeratorL}{0}
	\readarray{vNumeratorR}{0}

	\foreach \j in {1,...,#2}
	{        
		\pgfmathparse{ int(#4+\j+1) } \checkvKnots(\pgfmathresult)
		\pgfmathsetmacro{\numR}{\cachedata-#3}
		
		\pgfmathparse{ int(#4-\j+2) } \checkvKnots(\pgfmathresult)
		\pgfmathsetmacro{\numL}{#3-\cachedata} 					

		\expandarrayelementtrue
		\pgfmathparse{ int(\j+1) }
		\vNumeratorL(\pgfmathresult)={\numL}
		\vNumeratorR(\pgfmathresult)={\numR}
		
		\forloop[1]{itR}{0}{\value{itR} < \j }
		{
			\pgfmathparse{ int(\theitR+1) } \let\tS\pgfmathresult  	
			\checkvSave(\tS)							
			\pgfmathsetmacro{\save}{\cachedata}  
		
			\pgfmathparse{int(\j-\theitR+1)} \checkvNumeratorL(\pgfmathresult)
			\pgfmathsetmacro{\tmpL}{\cachedata} 	
		    
			\pgfmathparse{ int(\theitR+2) } \checkvNumeratorR(\pgfmathresult)
			\pgfmathsetmacro{\tmpR}{\cachedata} 	       
	
			\pgfmathparse{ int(\theitR+1) } \checkvN(\pgfmathresult)
			\pgfmathparse{ \cachedata / (\tmpL+\tmpR) } \let\tmp\pgfmathresult
			
			\pgfmathparse{ \save + \tmpR * \tmp } 
			\vN(\tS)={\pgfmathresult}

			\pgfmathparse{ \tmpL * \tmp } \let\tmpsave\pgfmathresult
			\pgfmathparse{ int(\theitR+2) } \let\tS\pgfmathresult      
			\vSave(\tS)={\tmpsave}
		}
		
		\pgfmathparse{ int(\j+1) } \checkvSave(\pgfmathresult )
		\vN(\tS)={\cachedata}
	}

	\pgfmathparse{int( #2+1) } \let\lastIndex\pgfmathresult
	\checkvN(1) \pgfmathsetmacro{\first}{\cachedata}  
	\foreach \i [remember=\a as \lasta (initially \first)] in {2,...,\lastIndex}
	{
		\checkvN(\i) \def\a{\lasta,\cachedata} 
		\ifthenelse{\i=\lastIndex}{ \xdef#5{\a} }{}
	}
    
	\foreach \i in {1,...,\lastIndex}
	{
	    \clrarray{vNumeratorR}(\i)
	    \clrarray{vNumeratorL}(\i)
	    \clrarray{vN}(\i)
	    \clrarray{vSave}(\i)
	}
	\delarray\vN
	\delarray\vNumeratorL
	\delarray\vNumeratorR
	\delarray\vSave
}

\newcounter{countvalues}
\newcounter{getBasis}
\newcommand{\bsplinebasis}[4] 
{						
    \newarray\vKnots 	
    \readarray{vKnots}{#1}

    \pgfmathparse{#3}  \let\i\pgfmathresult
    \pgfmathparse{#2}  \let\p\pgfmathresult

    \setcounter{countvalues}{0}
    \setcounter{getBasis}{\p} 
    \pgfmathparse{int(\i+\p)}
    \foreach \knotspan in {\i,...,\pgfmathresult}
    {
	\pgfmathparse{ int(\knotspan+1+1) } \checkvKnots(\pgfmathresult)
	\pgfmathsetmacro{\tmpR}{\cachedata} 						
	
	\pgfmathparse{ int(\knotspan+1) } \checkvKnots(\pgfmathresult) 
	\pgfmathsetmacro{\tmpL}{\cachedata} 						

	\pgfmathparse{ \tmpR - \tmpL } \let\spansize\pgfmathresult  			
   
        \pgfmathparse{ \spansize > 0.0 } \let\bNonZero\pgfmathresult
        \ifthenelse{ \bNonZero = 1 }
        {
            \foreach \percentU in {0,10,...,100}
            {
		\pgfmathparse{\tmpL+\spansize*\percentU/100} \let\u\pgfmathresult

		\bsplinevalue{#1}{\p}{\u}{\knotspan}{\Basis} 
		\def\basisfuncarray{{\Basis}} 				
		
		\pgfmathparse{\basisfuncarray[\thegetBasis]} 
		
		\ifthenelse{\thecountvalues=0}
		{ 
			\xdef\nodeB{"\u,\pgfmathresult"}
		}{
			\xdef\nodeB{\nodeB,"\u,\pgfmathresult"}  
		}
		\addtocounter{countvalues}{1}
            }
        }{}
        \addtocounter{getBasis}{-1}
    }
    
	\xdef#4{\nodeB}

	\delarray\vKnots
}

\newcommand{\bsplinebasisspan}[6] 
{							
    \newarray\vKnots 	
    \readarray{vKnots}{#1}

    \pgfmathparse{#5}  \let\plotknotspan\pgfmathresult
    \pgfmathparse{#4}  \let\splineknotspan\pgfmathresult
    \pgfmathparse{#3}  \let\i\pgfmathresult
    \pgfmathparse{#2}  \let\p\pgfmathresult

    \setcounter{countvalues}{0}
    \setcounter{getBasis}{\i} 
    \foreach \knotspan in {\plotknotspan}
    {
	\pgfmathparse{ int(\knotspan+1+1) } \checkvKnots(\pgfmathresult)
	\pgfmathsetmacro{\tmpR}{\cachedata} 						
	
	\pgfmathparse{ int(\knotspan+1) } \checkvKnots(\pgfmathresult) 
	\pgfmathsetmacro{\tmpL}{\cachedata} 						

	\pgfmathparse{ \tmpR - \tmpL } \let\spansize\pgfmathresult  			
   
        \pgfmathparse{ \spansize > 0.0 } \let\bNonZero\pgfmathresult
        \ifthenelse{ \bNonZero = 1 }
        {
            \foreach \percentU in {0,10,...,100}
            {
		\pgfmathparse{\tmpL+\spansize*\percentU/100} \let\u\pgfmathresult

		\bsplinevalue{#1}{\p}{\u}{\splineknotspan}{\Basis} 
		\def\basisfuncarray{{\Basis}} 				
		
		\pgfmathparse{\basisfuncarray[\thegetBasis]} 
		
		\ifthenelse{\thecountvalues=0}
		{ 
			\xdef\nodeB{"\u,\pgfmathresult"}
		}{
			\xdef\nodeB{\nodeB,"\u,\pgfmathresult"}  
		}
		\addtocounter{countvalues}{1}
            }
        }{}
        \addtocounter{getBasis}{-1}
    }
    
	\xdef#6{\nodeB}

	\delarray\vKnots
}

\newcommand{\plotbsplinebasis}[4] 	
{							
	\bsplinebasis{#1}{#2}{#3}{\nodeOut}
	\def\nodearray{{\nodeOut}}

	\xdef\name{ }
	\addtocounter{countvalues}{-1}
	\foreach \i in {0,...,\thecountvalues}
	{
		\pgfmathparse{\nodearray[\i]}
		\coordinate (point\i) at (\pgfmathresult);	  
		\xdef\name{ \name (point\i)  }
	}
	
	\draw[#4] plot coordinates{ \name };
	
	\xdef\name{ }
}

\newcommand{\plotbsplinesegment}[6] 	
{								
	\bsplinebasisspan{#1}{#2}{#3}{#4}{#5}{\nodeOut}
	\def\nodearray{{\nodeOut}}

	\xdef\name{ }
	\addtocounter{countvalues}{-1}
	\foreach \i in {0,...,\thecountvalues}
	{
		\pgfmathparse{\nodearray[\i]}
		\coordinate (point\i) at (\pgfmathresult);	  
		\xdef\name{ \name (point\i)  }
	}
	
	\draw[#6] plot coordinates{ \name };
	
	\xdef\name{ }
}

\newcommand{\plotbsplineaccumulated}[5] 	
{						
    \newarray\vKnots 	
    \readarray{vKnots}{#1}
    \newarray\vSubCoef 	
    \readarray{vSubCoef}{#3}
    
    \pgfmathparse{#4}  \let\plotknotspan\pgfmathresult
    \pgfmathparse{#4}  \let\splineknotspan\pgfmathresult
    \pgfmathparse{#2}  \let\p\pgfmathresult
    
    \setcounter{countvalues}{0}
    \pgfmathparse{int( \p+1) } \let\lastIndex\pgfmathresult
    \foreach \knotspan in {\plotknotspan}
    {
        \pgfmathparse{ int(\knotspan+1+1) } \checkvKnots(\pgfmathresult)
        \pgfmathsetmacro{\tmpR}{\cachedata} 						
        
        \pgfmathparse{ int(\knotspan+1) } \checkvKnots(\pgfmathresult) 
        \pgfmathsetmacro{\tmpL}{\cachedata} 						
        
        \pgfmathparse{ \tmpR - \tmpL } \let\spansize\pgfmathresult  			
        
        \pgfmathparse{ \spansize > 0.0 } \let\bNonZero\pgfmathresult
        \ifthenelse{ \bNonZero = 1 }
        {
            \foreach \percentU in {0,10,...,100}
            {
                \pgfmathparse{\tmpL+\spansize*\percentU/100} \let\u\pgfmathresult
                
                \bsplinevalue{#1}{\p}{\u}{\splineknotspan}{\Basis} 
                \def\basisfuncarray{{\Basis}} 				
                
                \setcounter{getBasis}{0} 
                \pgfmathparse{\basisfuncarray[\thegetBasis]} 
                \let\basisValue\pgfmathresult
                
                \checkvSubCoef(1) \pgfmathsetmacro{\coef}{\cachedata}  
                \pgfmathparse{ \basisValue * \coef } \let\first\pgfmathresult
                
                \xdef\lastx{\first}
                \foreach \i in {2,...,\lastIndex}
                { 
                    \addtocounter{getBasis}{1}       
                    \pgfmathparse{\basisfuncarray[\thegetBasis]}
                    \let\basisValue\pgfmathresult
                    
                    \checkvSubCoef(\i) \pgfmathsetmacro{\coef}{\cachedata}  
                    \pgfmathparse{ \lastx + \basisValue * \coef } \let\sum\pgfmathresult
                    
                    \xdef\lastx{\sum}
                    
                    \ifthenelse{\i=\lastIndex}
                    {
                        \ifthenelse{\thecountvalues=0}
                        { 
                            \xdef\nodeBB{"\u,\sum"}
                        }{
                            \xdef\nodeBB{\nodeBB,"\u,\sum"}  
                        }
                        \addtocounter{countvalues}{1}
                    }{}
                    
                }
            }
        }{}
    }
    
    \delarray\vKnots
    \delarray\vSubCoef
    
    \def\nodearray{{\nodeBB}}
    
    \xdef\name{ }
    \addtocounter{countvalues}{-1}
    \foreach \i in {0,...,\thecountvalues}
    {
        \pgfmathparse{\nodearray[\i]}
        \coordinate (point\i) at (\pgfmathresult);                
        \xdef\name{ \name (point\i)  }
    }
    
    \draw[#5] plot coordinates{ \name };
    
    \xdef\name{ }
}

\begin{document}
    
\title{Efficient isogeometric Boundary Element simulation of elastic domains containing thin inclusions}
\begin{frontmatter}

\author[unife]{Vincenzo Mallardo\corref{cor1}}
\author[ifbaddr]{Christian Dünser}
\author[ifbaddr]{Gernot Beer}

\address[unife]{Department of Architecture, University of Ferrara, Via Quartieri 8, 44121 Ferrara, Italy}

\address[ifbaddr]{Institute of Structural Analysis, Graz University
  of Technology, Lessingstraße 25/II, 8010 Graz, Austria}

\cortext[cor1]{Corresponding author.
  Tel.Fax: +39 0532 293621, mail: \url{mlv@unife.it}, web: \url{http://docente.unife.it/docenti-en/vincenzo.mallardo?set_language=en}}

\begin{abstract}
This paper is concerned with the Boundary Element simulation of elastic domains that contain thin inclusions that have elastic material properties, which are different to the domain. 
With thin inclusions we mean inclusions with extreme aspect ratios, i.e. where one dimension is much smaller than the other ones.
Examples of this are reinforcements in civil/mechanical engineering or concrete linings in underground construction. 
The fact that an inclusion has an extreme aspect ratio poses a challenge to the numerical integration of the arising singular integrals and novel approaches are presented to deal with it. Several examples demonstrate the efficiency and accuracy of the proposed methods and show that the results are in good agreement with analytical and other numerical solutions.
\end{abstract}

\end{frontmatter}

\section{Introduction}
Isogeometric methods for simulation were first introduced in \cite{Hughes2005a} and were well received in the Finite Element Method (FEM) community. First implementations of the isogeometric Boundary Element Method (BEM) followed later, first in 2-D \cite{simpson2012two} and then in 3-D \cite{scott2013isogeometric} elasticity. A book on this topic was published recently \cite{BeerMarussig}. It was found that isogemetric methods lead to more accurate, user friendly and efficient simulations.

Because of its boundary-only modelling, the combination BEM – NURBS (Non Uniform Radial Basis Functions) is an ideal choice for isogeometric analysis (IgA) of solids. Many applications are provided in literature. 

For instance in \cite{Han_101} the isogeometric analysis in boundary element method (IGABEM) is tested for 2D potential problems. As the computation of singular integrals is a key issue, \cite{Han_101} develops a semianalytical scheme to determine the nearly singular integrals that occur with thin problems (such as thin-coating) and with the measurement of physical quantities of near boundary points. A new adaptive IGABEM with hierarchical B-splines is proposed in \cite{Falini_102}. A local quadrature scheme to solve 2D Laplace problems is developed to compute nearly singular integrals. The IGABEM can be accelerated by using suitable Fast Multipole Method (FMM) schemes. In \cite{Li_103} for instance the black box FMM is adopted for the solution of 3D elastostatic analysis.
In shape optimization the BEM represents a natural choice: it minimizes the mesh generation/regeneration burden and it maximizes the benefit to directly use an existing CAD geometry. In \cite{Lian_104} a gradient based shape optimization appraoch is developed in 2D linear, homogeneous, isotropic materials.
The advantages of FEM and BEM are tested in \cite{Liu_105} by adopting an isogeometric approach for structural-acoustic analysis problems. The infinite fluid domain is modelled through a BEM discretization whereas the shell is represented by FEM. The isogeometric approach allows the direct use of high order discretisations generated through CAD software without any effort in model conversion.

It must be pointed out that many elasticity problems deal with inhomogeneous domains. On the other hand, the basic BEM can only deal with homogeneous domains and this is a severe restriction to the application of the method to real world problems. Several publications have addressed this problem and proposed solutions. 

A very basic formulation in dealing with this is to discretise the boundary of the matrix and the surface of each inclusion and to impose continuity and equilibrium conditions on the interfaces (multi-region BEM \cite{Banerjee1981,Achenbach1990}). The mechanics of the unilateral contact can be also included by, for instance, a linear complementarity approach \cite{Mallardo2000}. 
Many drawbacks arise with such an approach, mainly 1) the governing system of equations grows up hugely with increasing inclusions, 2) numerical issues (singularity) occur with thin inclusions and 3)
the direct coupling may be too computationally demanding for large models. An adaptive iterative coupling algorithm is proposed in \cite{Wang_121} in the context of the Laplace equation to alleviate such an issue, but still many iterations are required, and no testing is carried out for elasticity.

An alternative approach consists in dealing with non linear and/or non-homogeneous problems by supplementing the boundary integral equation with a domain integral. A comparison of domain integral evaluation techniques is provided in \cite{Ingber_110} with reference to the Poisson equation. Initially, some methods were proposed to avoid the domain discretization, notably dual and multiple reciprocity methods \cite{partridge1992} and particular solution methods \cite{Banerjee1981}. All these methods require both approximating a function in the interior of the domain and a corresponding particular solution, but they suffer in precision and efficiency as being sensitive to the location of the radial basis functions, showing possible ill-behavior of the particular solution and being problem-dependent.
A novel approach is proposed in \cite{Gao_122}. The multi-medium elastic problems are solved by adding an interface integral and a whole domain integral to the basic integral equation, therefore internal cells in the whole domain and boundary elements on the interfaces are required along with the classical external boundary discretization. Furthermore, thin inclusions introduce additional singularity problems.

The mechanical analysis of multi-medium domains is therefore an open issue: the computational time associated to the geometry generation and to the mesh discretization is still unsuitable for large scale real problems and the presence of thin inclusions, with related singularity issues, have not been investigated so far.
The authors have published previously on this topic (see for example \cite{Beer2016} and \cite{Beer17}), introducing efficient and user friendly ways for dealing with the volume integrals, including avoiding the additional discretisation effort via cells. The problem with the methods published there is that they do not work for thin inclusions, i.e. for inclusions where one of the dimensions is much smaller than the other dimensions. 
The reason for this is that the numerical integration of the singular Kernel functions poses problems. Special numerical techniques are required and this is the main topic of the paper.

We first introduce the governing integral equations and how they can be solved numerically. For a better understanding of the paper we repeat some already published material here. 
The numerical integration of the volume integral is discussed in detail.
First we show how the inclusion geometry can be defined with few parameters using NURBS basis functions.
To make the simulation of thin inclusions efficient and accurate we introduce a local coordinate system where the strains and stresses are defined.
We also show how displacements can be interpolated using NURBS and how derivatives can be obtained.
The efficient numerical integration of the singular functions is addressed. Two methods are presented: One where the integration is purely numerical and one where it is a mixture of numerical and analytical evaluation. The latter is required for very thin inclusions, where the purely numerical approach breaks down.

It is demonstrated on three test examples that the proposed implementation works. In the first simple we test the limit of the purely numerical integration. In the second we show that good results can also be obtained for more complex problems. Finally, in the last example the method is applied to a curved geometry to demonstrate the efficiency introduced by the introduction of local strains.

\section{Theory}
To account for the influence of inclusions that have elastic properties that are different to the domain, we use the method of initial stress, well known in the FEM community for dealing with plasticity.
Consider the elastic stresses $\myVecGreek{\sigma}$ due to strains $\myVecGreek{\epsilon} $:
\begin{equation}
\label{ }
 \myVecGreek{\sigma}= \mathbf{ D} \myVecGreek{\epsilon} 
\end{equation} 
where $ \mathbf{ D}$ is the elasticity matrix.

The initial stresses $ \myVecGreek{\sigma}_{0} $ due to a difference in elastic properties between the domain and the inclusions are given by:
\begin{equation}
\label{InitialS}
\myVecGreek{\sigma}_{0} =( \mathbf{ D} - \mathbf{ D}_{incl}) \myVecGreek{\epsilon} =\mathbf{ D}^{\prime}\myVecGreek{\epsilon} 
\end{equation}
where $ \mathbf{ D}_{incl}$ is the elasticity matrix for the inclusion.

\subsection{Governing integral equations}
Consider a domain  $\domain$ bounded by $\boundary$ with a subdomain $\domain_{0}$, where initial stresses $ \myVecGreek{\sigma}_{0}(\fieldpt)$ are present.

To establish the integral equations we apply Betti's theorem and the collocation method.
The regularised integral equations are written in matrix notation as:
\begin{equation}
\label{Inteq}
    \begin{aligned}
\int_{\boundary} \fund{T}(\sourcept_{n},\fieldpt) ( \myVec{u}(\fieldpt) - \myVec{u}(\sourcept_{n})) \ d\boundary(\fieldpt) - \mathbf{ A}_{n} \myVec{u}(\sourcept_{n}) &=& \int_{\boundary} \fund{U}(\sourcept_{n},\fieldpt) \ \myVec{t}(\fieldpt) \ d\boundary(\fieldpt) \\
+  \int_{\domain_{0}} \fund{E} (\sourcept_{n},\fieldpt)  \myVecGreek{\sigma}_{0} (\fieldpt) d \domain_{0} (\fieldpt)  .
\end{aligned}
\end{equation}
where $\fund{T}(\sourcept_{n},\fieldpt)$ is a matrix that contains the fundamental solutions for the tractions, $\fund{U}(\sourcept_{n},\fieldpt)$ the one for the displacements and $\fund{E} (\sourcept_{n},\fieldpt)$ for the strains at point $\fieldpt$ due to a source at $\sourcept_n$.
$\myVec{u}$ and $\myVec{t}$ are displacement and traction vectors at boundary points.
 $\mathbf{ A}_{n}$ represents the azimuthal integral that is equal to the unit matrix for infinite and the zero matrix for finite domain problems.
The fundamental solutions for the displacements and tractions are well published (see for example \cite{BeerMarussig}).

The fundamental solution $\fund{E}$ is in index notation:
\begin{equation}
\label{ }
E_{ijk}= \frac{-C}{r^{2}}\left[C_{3}(r_{,k}\delta_{ij} + r_{,j}\delta_{ik}) - r_{,i}\delta_{jk} + C_{4} \ r_{,i} r_{,j} r_{,k}\right]
\end{equation}
where the index $i$ specifies the direction of the source, $jk$ the strain component. The constants are: $C=\frac{1}{16 \pi G (1- \nu})$, $C_{3}=1 - 2\nu$ and $C_{4}=3$. 

The equation can be expressed as:
\begin{equation}
\label{ }
E_{ijk}= \frac{1}{r^{2}}\tilde{E}_{ijk}
\end{equation}
where
\begin{equation}
\label{ }
\tilde{E}_{ijk}= -C\left[C_{3}(r_{,k}\delta_{ij} + r_{,j}\delta_{ik}) - r_{,i}\delta_{jk} + C_{4} \ r_{,i} r_{,j} r_{,k}\right]
\end{equation}

Using Voight notation for the initial stresses 
\begin{equation}
\label{Voightnot}
\myVecGreek{\sigma}_{0}= \left \{\begin{array}{c}\sigma_{x0} \\\sigma_{y0} \\\sigma_{z0} \\\tau_{xy0} \\\tau_{yz0}\\\tau_{xz0}\end{array}\right\}
\end{equation}
we convert the tensor $E_{ijk}$ to a matrix $\fund{E}$:
\begin{equation}
\label{eq_Voigt_convertion}
\fund{E}= \left[\begin{array}{cccccc}E_{111} & E_{122} & E_{133}  & E_{112}  + E_{121}& E_{123} + E_{132} & E_{113} +E_{131}  \\E_{211} & E_{222} & E_{233}  & E_{212} + E_{221}  & E_{223}  + E_{232} & E_{213} + E_{231}\\E_{311} & E_{322} & E_{333}  & E_{312} + E_{321} & E_{323}  + E_{332}& E_{313} + E_{331}\end{array}\right]
\end{equation}
The integral equations \ref{Inteq} are solved numerically using Gauss Quadrature. For this a discretisation of the surface and volume integrals is necessary.

The discretisation involves 2 steps:
\begin{itemize}
  \item  The subdivision of the boundary domain into patches and the volume domain into inclusions
  \item  The approximation of the unknown boundary values and the approximation of initial stresses.
\end{itemize}
This will be discussed in the subsequent sections separately for the boundary and volume integrals.

\section{Discretisation of the boundary integrals}
For the numerical solution of the boundary integral equations the integrals are expressed as sum of integrals over patches:
\begin{equation}
    \begin{aligned}
        \label{Regular2}
 \int_{\boundary} \fund{U}(\sourcept_{n},\fieldpt) \ \myVec{t}(\fieldpt) \ d\boundary(\fieldpt) - \int_{\boundary} \fund{T}(\sourcept_{n},\fieldpt) ( \myVec{u}(\fieldpt) - \myVec{u}(\sourcept_{n})) \ d\boundary(\fieldpt) + \mathbf{ A}_{n} \myVec{u}(\sourcept_{n})  \approx  \\
        \sum_{e=1}^{E} \int_{\Gamma_{e}}  \fund{U}(\sourcept_{n},\fieldpt) \ \myVec{\dual}^{e}(\fieldpt)\ d\Gamma_{e}(\fieldpt)  
        - \sum_{e=1}^{E} \int_{\Gamma_{e}} \fund{T}(\sourcept_{n},\fieldpt)  \myVec{\primary}^{e}(\fieldpt) d \Gamma_{e}  \\
        +  \left[\sum_{e=1}^{E} \ \left(\int_{\Gamma_{e}} \fund{T}(\sourcept_{n},\fieldpt) d \Gamma_{e} \right) + \mathbf{A}_{n}\right] \myVec{\primary}(\sourcept_{n})
    \end{aligned}
\end{equation}
where $e$ specifies the patch number and $E$ is the total number of patches. In the following the geometry of patches is specified using NURBS basis functions. The advantage of this is that some geometrical shapes such as cylinder and spheres can be described exactly with few parameters.
For further information on NURBS and how (\ref{Regular2}) is obtained the reader is referred to \cite{BeerMarussig}

There are 2 types of patches that we use in this paper: finite and infinite patches. 

\subsection{Geometry definition of finite patches}
\begin{figure}
\begin{center}
\begin{overpic}[scale=0.5]{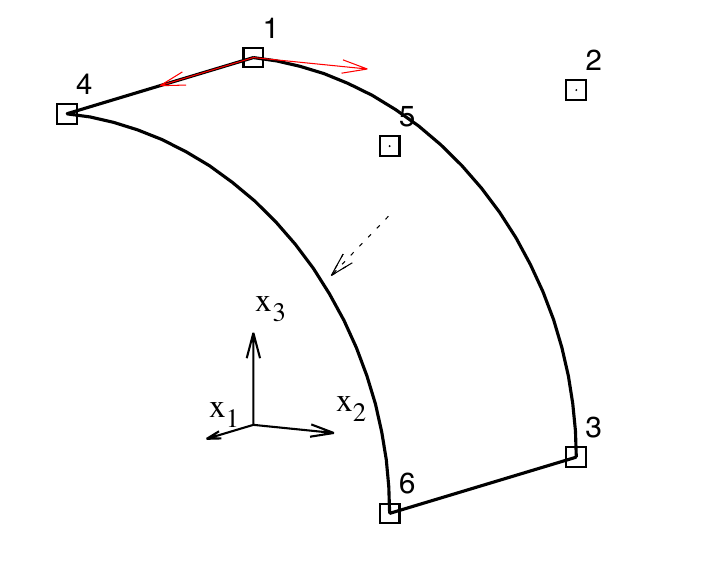}
 \put(50,75){$\uu$}
 \put(20,70){$\vv$}
\end{overpic}
\begin{overpic}[scale=0.5]{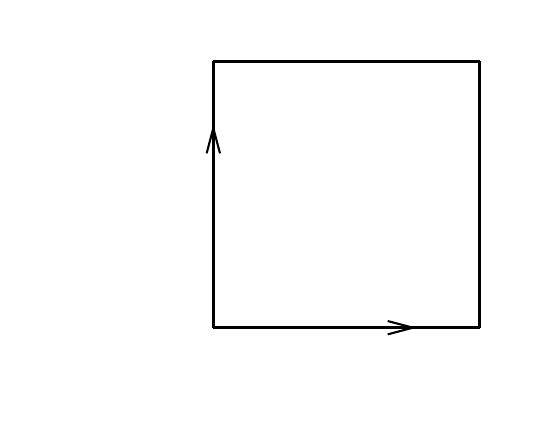}
 \put(70,25){$\uu$}
 \put(30,50){$\vv$}
\end{overpic}
\caption{A finite patch with control points (numbered squares). Left: in the global, right: in the local coordinate system. Also shown is the ``outward normal''.}
\label{fig6:Patch3D}
\end{center}
\end{figure}

In \myfigref{fig6:Patch3D} we show an example of a finite patch.
The mapping from the local $\myVecGreek{\xi}(\uu,\vv)$ to the global $\pt{x}$ coordinate system is given by
\begin{equation}
\pt{x} (\xi,\eta)= \sum_{i=1}^{I} \NURBS_{i} (\xi, \eta) \pt{x}_{i}.
\end{equation}
where $\NURBS_{i} (\xi, \eta)$ are NURBS basis functions and the control points (coordinates $\pt{x}_{i}$) are numbered consecutively, first in the $\uu$- and then in the $\vv$-direction.

The vectors tangential to the surface are given by

\begin{align}
\mathbf{ V}_{\xi}= \frac{\partial \pt{x}}{\partial \xi}= \left(\begin{array}{c}\frac{\partial x}{\partial \xi} \\ \\ \frac{\partial y}{\partial \xi} \\ \\ \frac{\partial z}{\partial \xi}\end{array}\right)
&& \text{and} && \mathbf{ V}_{\eta}= \frac{\partial \pt{x}}{\partial \eta}= \left(\begin{array}{c}\frac{\partial x}{\partial \eta} \\ \\ \frac{\partial y}{\partial \eta} \\ \\ \frac{\partial z}{\partial \eta} \end{array}\right)
\end{align}
and the unit vector normal is
\begin{equation}
\mathbf{ n}= \frac{\mathbf{V}_{\xi} \times  \mathbf{ V}_{\eta}}{J}.
\end{equation}
The Jacobian is
\begin{equation}
J= | \mathbf{ V}_{\xi} \times  \mathbf{ V}_{\eta} |.
\end{equation}
 The  direction of the ``outward normal'' depends on how the control points are numbered.

\subsection{Geometry definition of infinite patches}
\label{sec:InfinitePatches}
\begin{figure}
\begin{center}
\begin{overpic}[scale=0.5]{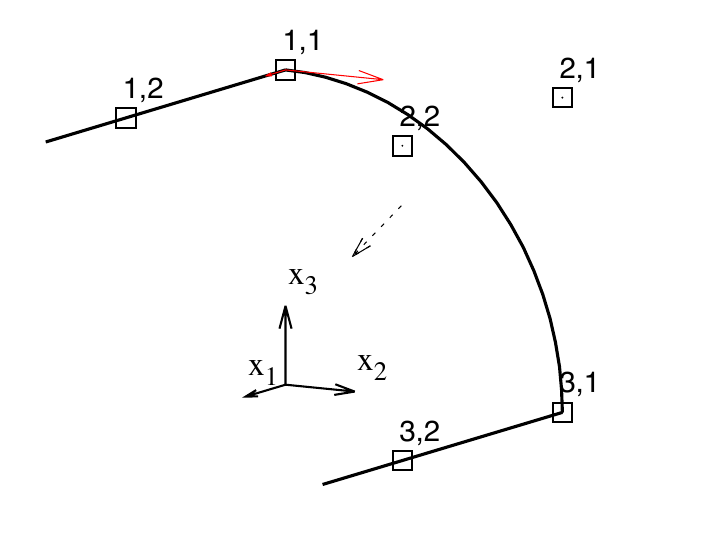}
 \put(50,70){$\uu$}
 \put(30,68){$\vv$}
\end{overpic} 
\begin{overpic}[scale=0.5]{pics/Patch3Dst.pdf}
 \put(70,25){$\uu$}
 \put(30,50){$\vv$}
\end{overpic}
\caption{Example of an infinite patch. Left in the global and right in the local coordinate system}
\label{fig6:Ptch3Dinf}
\end{center}
\end{figure}
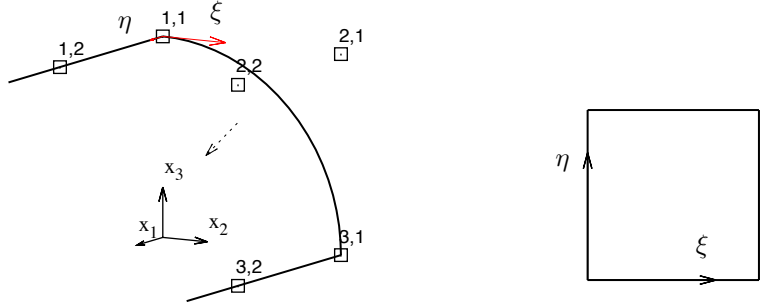

Here we introduce a patch definition that is useful for the simulation in geomechanics where one sometimes has to consider a surface that tends to infinity \cite{Beer2015b}. In this case we define an infinite patch as shown in \myfigref{fig6:Ptch3Dinf}.
The mapping for a patch that extends to infinity in the $\eta$-direction is given by
\begin{equation}
\pt{x}= \sum_{j=1}^{2}\sum_{i=1}^{I} \NURBS^{\infty}_{ij}(\xi,\eta) \pt{x}_{ij}
\end{equation}
where

\begin{equation}
\NURBS^{\infty}_{ij}(\xi,\eta)= \NURBS_{i}(\xi) M_{j}^{\infty}(\eta)
\end{equation}

and the special infinite basis functions are

\begin{align}
 M_{1}^{\infty} =  \frac{1 - 2\eta}{1-\eta} && \text{and} &&  M_{2}^{\infty}  =  \frac{\eta}{1-\eta}. 
\end{align}

The vectors in the tangential directions are given by

\begin{eqnarray}
\mathbf{V}_{\xi}= \frac{\partial \pt{x}}{\partial \xi}=\sum_{j=1}^{2}\sum_{i=1}^{I} \frac{\partial \NURBS_{i}(\xi)}{\partial \xi}  M_{j}^{\infty}(\eta) \pt{x}_{ij} \\
\mathbf{ V}_{\eta}= \frac{\partial \pt{x}}{\partial \eta}=\sum_{j=1}^{2}\sum_{i=1}^{I} \NURBS_{i}(\xi)  \frac{\partial M_{j}^{\infty}(\eta)}{\partial \eta} \pt{x}_{ij} 
\end{eqnarray}

where

\begin{align}
\frac{\partial M_{1}^{\infty}}{\partial \eta}= \frac{-1}{(1-\eta)^{2}} && \text{and} && \frac{\partial M_{2}^{\infty}}{\partial \eta}= \frac{1}{(1-\eta)^{2}}.
\end{align}

The unit vector normal is computed as for the finite patch. It is noted that the Jacobian $J$ tends to infinity as $\vv$ tends to 1.

\newpage

\section{Geometrical discretisation of Volume integrals}
\label{Geom}
\begin{figure}
\begin{center}
\begin{overpic}[scale=1]{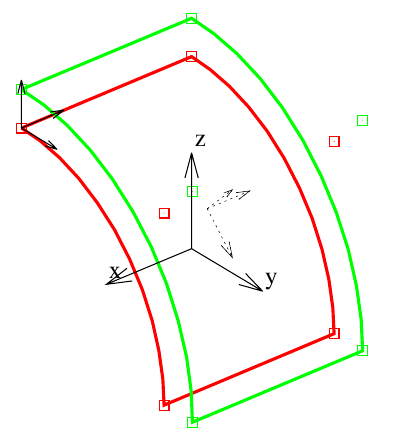} 
\put(60,55){$y^{\prime}$}
\put(55,42){$x^{\prime}$}
\put(50,60){$z^{\prime}$}
\put(12,68){$s$}
\put(15,75){$t$}
\put(5,85){$r$}
\end{overpic}
\begin{overpic}[scale=0.5]{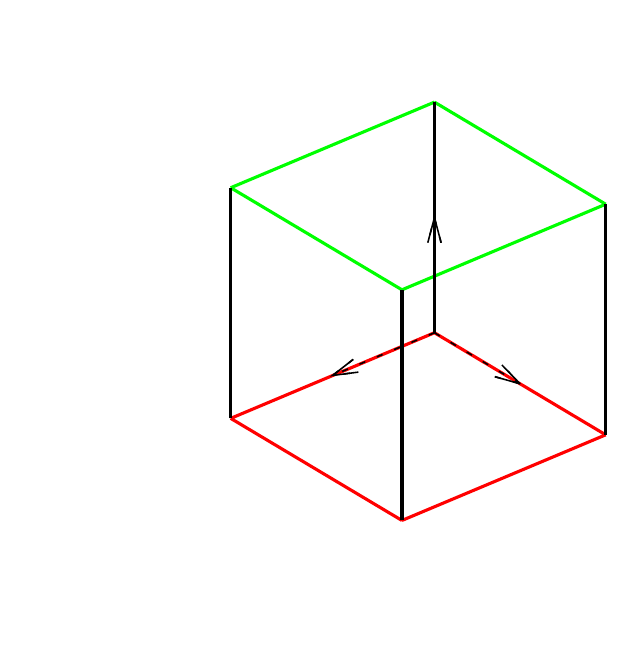} 
\put(50,50){$s$}
\put(75,50){$t$}
\put(60,70){$r$}
\end{overpic}
\caption{An inclusion showing coloured bounding surfaces. Right: in local $\pt{s}$ coordinate system. Left: in global $\pt{x}$ coordinate system with control points and local axes.}
\label{Shotcrete}
\end{center}
\end{figure}
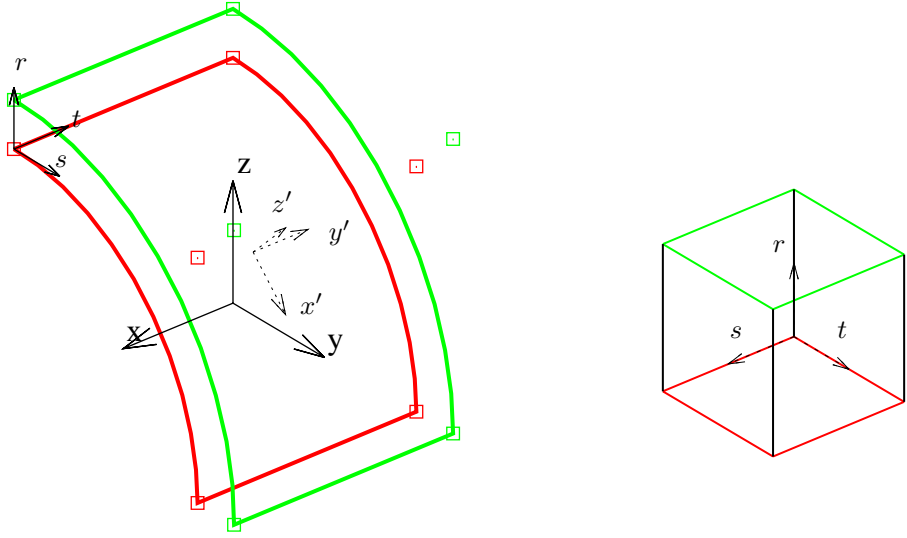

For the geometrical discretisation of the volume integral we subdivide $\domain_0$ into one or several \textbf{inclusions}.
The geometry of an inclusion is described by 2 bounding NURBS surfaces. 
In \myfigref{Shotcrete} we show the colour coded definition of the bounding surfaces as well as the associated control points.

\subsection{Mapping}
The global coordinates of a point $\pt{x}$ with the local coordinates $\pt{s}=(s,t,r)^{\mathrm{T}}=[0,1]^3$ are given by
\begin{equation}
\label{Mapx}
\pt{ x}({s,t,r})= (1-r) \ \pt{ x}^{I}(s,t) + {r} \ \pt{ x}^{II}({s,t})
\end{equation}
where
\begin{align}
  \pt{ x}^{I}({s,t})=\sum_{k=1}^{K} R_{k}({s,t}) \ \pt{ x}_{k}^{I} && \text{and} &&\pt{ x}^{II}({s,t})=\sum_{k=1}^{K} R_{k}({s,t}) \ \pt{ x}_{k}^{II} .
\end{align}
The superscript $I$ relates to the bottom (red) surface and $II$ to the top (green) surface and $ \pt{ x}_{k}^{I} $, $ \pt{ x}_{k}^{II} $ are control point coordinates. $K$  represents the number of control points which is the same for each surface, $R_{k}({s,t})$  are NURBS basis functions. 

The vectors in local directions are given by:
\begin{equation}
\label{Vec}
\begin{aligned}
 \mathbf{ V}_{s} = \frac{\partial \pt{ x}({s,t,r})}{\partial {s}}&=& (1-{r}) \ \frac{\partial \pt{ x}^{I}({s,t})}{\partial {s}} &+&& {r} \ \frac{\partial \pt{ x}^{II}({s,t}) }{\partial {s}} \\  \mathbf{ V}_{t} =   \frac{\partial \pt{ x}({s,t,r})}{\partial {t}}&=& (1-{r}) \ \frac{\partial \pt{ x}^{I}({s,t})}{\partial {t}} &+&& {r} \ \frac{\partial \pt{ x}^{II}({s,t}) }{\partial {t}} \\
 \mathbf{ V}_{r}= \frac{\partial \pt{ x}({s,t,r})}{\partial {r}}&=& -\pt{ x}^{I}({s,t}) &+&&
  \ \pt{ x}^{II}({s,t})
\end{aligned}
\end{equation}
where
\begin{align}
 \frac{\partial\pt{ x}^{I}({s,t})}{\partial {s}}=\sum_{k=1}^{K} \frac{\partial R_{k}({s,t})}{\partial {s}} \ \pt{ x}_{k}^{I} &&\text{and} &&
 \frac{\partial\pt{ x}^{II}({s,t})}{\partial {}s}=\sum_{k=1}^{K} 
  \frac{\partial R_{k}({s,t})}{\partial {s}} \ \pt{ x}_{k}^{II}  \\
  \frac{\partial\pt{ x}^{I}({s,t})}{\partial {t}}=\sum_{k=1}^{K} \frac{\partial R_{k}({s,t})}{\partial {t}} \ \pt{ x}_{k}^{I} &&\text{and} &&
 \frac{\partial\pt{ x}^{II}({s,t})}{\partial {t}}=\sum_{k=1}^{K} 
  \frac{\partial R_{k}({s,t})}{\partial {t}} \ \pt{ x}_{k}^{II}
\end{align}

The Jacobi matrix of this mapping is
\begin{equation}
\label{ }
\mathbf{ J}= \left[\begin{array}{c} \mathbf{ V}_{s} \\ \mathbf{ V}_{t} \\ \mathbf{ V}_{r}\end{array}\right]
\end{equation}
and the Jacobian is $J(\pt{s})=| \mathbf{ J} |$.

\

\remark{For thin inclusions, where the thickness ($d$) is constant, only one bounding surface is required to describe the geometry. In this case we set $r=0$ in \myeqref{Mapx} and (\ref{Vec}). The numerical integration is carried out over the surface and multiplied by $d$. In this case $\mathbf{ V}_{r}= \mathbf{ n} \ d$ where $\mathbf{ n}$ is the unit vector normal to the surface.}

\

\subsection{Local coordinate system for defining stresses/strains}
\label{Loco}
Vectors $\mathbf{ V}_{s}$ and  $\mathbf{ V}_{t}$ may not be orthogonal.
In order to define local stresses and strains we require an orthogonal system.
We define a local coordinate system whereby the $x^{\prime}$ and $y^{\prime}$ axes are orthogonal and tangential to the surface and $z^{\prime}$ is normal to it. 

We assume the vector in $x^{\prime}$-direction to be in the direction of $s$:
\begin{equation}
\label{ }
\mathbf{ V}_{1}= \mathbf{ V}_{s}
\end{equation}
The vector in $z^{\prime}$-direction is given by:
\begin{equation}
\label{ }
\mathbf{ V}_{3}= \mathbf{ V}_{s} \times  \mathbf{ V}_{t}
\end{equation}
and the vector in $y^{\prime}$-direction is then obtained
\begin{equation}
\label{ }
\mathbf{ V}_{2}= \mathbf{ V}_{1} \times \mathbf{ V}_{3}
\end{equation}

The unit vectors in the local coordinate directions are:
\begin{eqnarray}
\mathbf{ v}_{1}= \frac{\mathbf{ V}_{1}}{|\mathbf{ V}_{1}|}  \text{   ,   } \mathbf{ v}_{2}= \frac{\mathbf{ V}_{2}}{|\mathbf{ V}_{2}|}   \text{   ,   } \mathbf{ v}_{3}= \frac{\mathbf{ V}_{3}}{|\mathbf{ V}_{3}|} 
\end{eqnarray}

\newpage

\section{Approximation of the unknown boundary values}
\begin{figure}
\begin{center}
\begin{overpic}[scale=0.5]{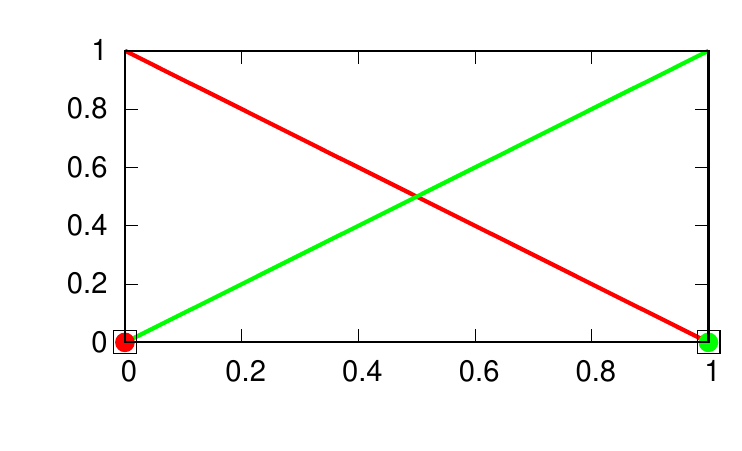}
\put(0,30){$R_{i}$}
\put(50,2){$\xi$}
\end{overpic}
\begin{overpic}[scale=0.5]{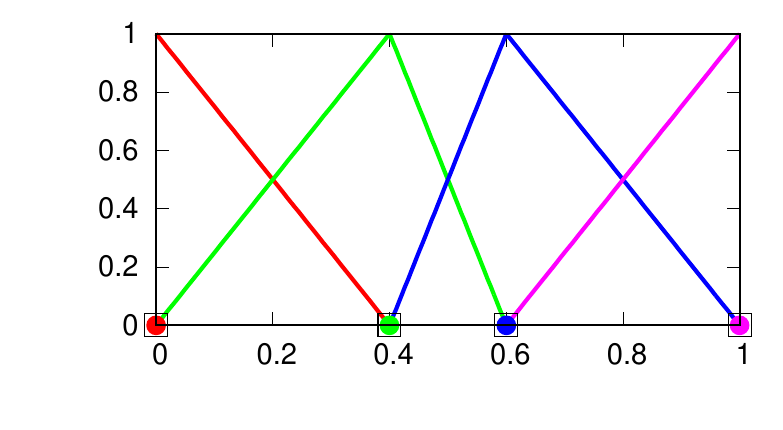}
\put(3,30){$\hat{R}_{i}$}
\put(60,0){$\xi$}
\end{overpic}
\caption{Example of generation of basis functions for approximation of the unknown. Left: Linear basis functions for the geometry, Right: Refined functions by 2 knot insertions at 0.4 and 0.6}
\label{o1k2}
\end{center}
\end{figure}

\begin{figure}
\begin{center}
\begin{overpic}[scale=0.5]{pics/basiso1.pdf}
\put(0,30){$R_{i}$}
\put(50,2){$\xi$}
\end{overpic}
\begin{overpic}[scale=0.5]{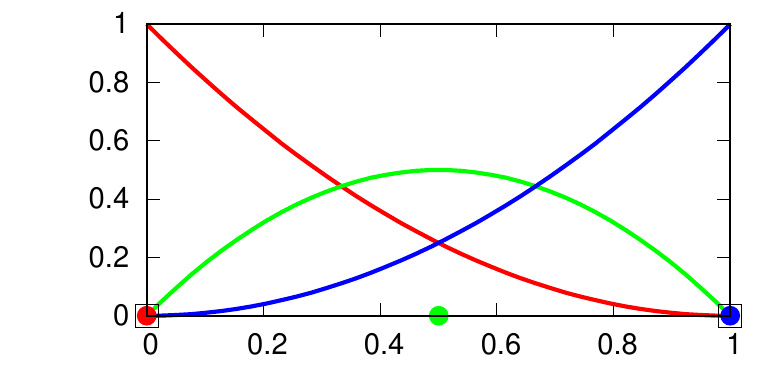}
\put(3,30){$\hat{R}_{i}$}
\put(60,0){$\xi$}
\end{overpic}
\caption{Example of generation of basis functions for approximation of the unknown. Left: Linear basis functions for the geometry, Right: Refined functions by one order elevation}
\label{o1o2}
\end{center}
\end{figure}
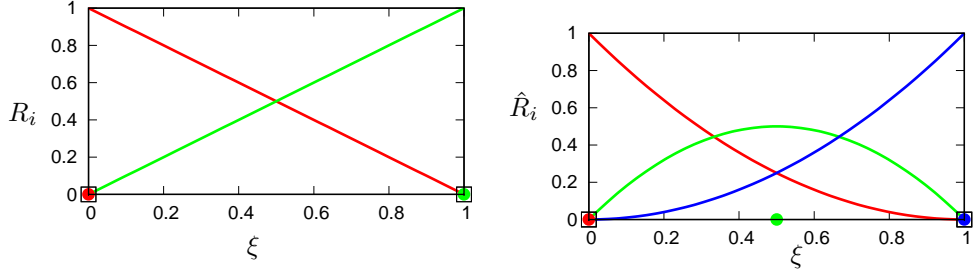

We approximate the unknown boundary values also with NURBS basis functions. For example for the displacements we have:
\begin{equation}
\mathbf{ u}(\xi,\eta)= \sum_{i=1}^{I}\hat{ \NURBS_{i} }(\xi, \eta) \mathbf{ u}_{i}.
\end{equation}
where $\hat{\NURBS}_{i} (\xi, \eta)$ are NURBS basis functions and $\mathbf{ u}_{i}$ are parameter values. 
We use the \textit{independent field approximation method} already published and verified in \cite{doi:10.1002/nme.5778,Marussig2015}, i.e. we leave the definition of the geometry untouched and only enrich the basis functions for the approximation of the unknown. This means that $\hat{ \NURBS_{i} }(\xi, \eta)$ are basis functions that are obtained by refining the basis functions $\NURBS_{i} (\xi, \eta)$ for the geometry.
Examples of the refinement of basis function space are shown in \myfigref{o1k2} and \myfigref{o1o2}.
In the first case two knots are inserted and in the second case the order is elevated by one. Also shown are the anchors of the basis functions which are used to determine the collocation points $\sourcept_n$.

For the infinite patches we assume plane strain conditions along the direction to infinity, i.e. we assume that the displacements are constant in this direction.

\section{Approximation of values inside an inclusion}
For the approximation of values inside an inclusion we establish a regular grid inside it and store values of displacement  in vectors $\{\mathbf{  u} \}$ and values of initial stress in $\{ \myVecGreek{\sigma}_{0}\} $ at all grid points.
The value of initial stress at any point with the local coordinates $\pt{s}$ is obtained by interpolation between grid points.
\begin{equation}
\label{eq11:Interpol}
 \myVecGreek{\sigma}_{0}(\pt{s})= \sum_{j=1}^{J} M_{j}^{\sigma}(\pt{s})  \myVecGreek{\sigma}_{0j}
\end{equation}
where $  \myVecGreek{\sigma}_{0j} $ is the initial stress at grid point $j$ with the local coordinate $\pt{s}_{j}$. $M_{j}^{\sigma}(\pt{s})$ are piecewise constant or linear basis functions. 
For the approximation of $\mathbf{ u}$ we use Lagrange polynomials or NURBS as will be explained later.

\section{Discretised integral equations}

Introducing shape functions for the variation of the unknown into \myeqref{Regular2}, the following discrete system of equations can be obtained \cite{BeerMarussig}:
\begin{equation}
\label{DisIE}
[\mathbf{ L}] \{\mathbf{ x}\} = \{\mathbf{ r}\} + [\mathbf{ B}_{0}] \{ \myVecGreek{\sigma}_{0}\} 
\end{equation}
where  $[\mathbf{ L}] $ is the left hand side, $\{\mathbf{ x}\}$ is the vector of unknowns and $ \{\mathbf{ r}\}$ is the right hand side for the system without inclusion . 

The sub-matrices of matrix $ [\mathbf{ B}_{0}]$, related to collocation point $n$ and grid point $j$, are given by
\begin{equation}
\label{ }
\mathbf{ B}_{0nj} = \sum_{m=1}^{M} \ \int_{\domain_{m}}  \fund{E} (\sourcept_{n},\fieldpt) M_{j}^{\sigma} (\fieldpt) d \domain_{m} (\fieldpt) 
\end{equation}
where $M$ is the total number of inclusions, $\domain_{m}$ is the volume of an inclusion and $ M_{j}^{\sigma} $ are interpolation functions introduced earlier.

For the evaluation of the surface integrals we use established and well published numerical procedures based on Gauss Quadrature (see for example \cite{BeerMarussig}). The evaluation of the volume integrals will be discussed later.

\section{Computation of values at grid points inside the inclusion}
To obtain the initial stresses we need to compute the displacements and strains at points inside the inclusion.
The displacement $ \mathbf{ u} $ at a point $\pt{x}$ inside the inclusion is given by:
\begin{equation}
    \begin{aligned}
        \label{eps0}
       \mathbf{ u}(\pt{x}) &= \int_{\boundary} \left[ \fund{U}(\pt{x},\fieldpt) \   \myVec{\dual} (\fieldpt)  - \fund{T}(\pt{x},\fieldpt) \   \myVec{\primary} (\fieldpt) \right] d\boundary (\fieldpt) \\
        &+ \int_{\domain_{0}} \fund{E} (\pt{x},\fieldpt) \myVecGreek{\sigma}_{0} (\fieldpt)  d \domain_{0} (\fieldpt)   
    \end{aligned}
\end{equation}

Displacement vectors at all grid points are now gathered in $\{\mathbf{ u} \}$. In matrix form we have:
\begin{equation}
\label{ }
\{\myVec{u}\}= [\hat{\mathbf{ A}}] \{\mathbf{ x}\} + \{\bar{\mathbf{ c}}\} + [\bar{\mathbf{ B}}_{0}]\{ \myVecGreek{\sigma}_{0}\}   
\end{equation}
where $[\hat{\mathbf{ A}}]$ is an assembled matrix that multiplies with the unknown $\{\mathbf{ x}\}$ and $\{\bar{\mathbf{ c}}\}$ collects the displacements due to given BC's. $ [\bar{\mathbf{ B}}_{0}]$ is similar to $ [\mathbf{ B}_{0}]$ except that the grid point coordinates $\pt{x}_{i}$ replace the source point coordinates $\sourcept_{n}$.

Because of the singularity of $\fund{T}$ the displacements can not be computed on the boundary. For points on a patch boundary ($\pt{x}_{k}$) we use:
\begin{equation}
\label{ }
\mathbf{ u}(\pt{x}_{k})=\sum_{i}^{I} \hat{R}_{i}(\myVecGreek{\xi}_{k}) \mathbf{ u}_{i}^{e}
\end{equation}
where $\hat{R}_{i}(\myVecGreek{\xi})$ are the NURBS basis functions used for approximating the displacements in patch $e$. The superscript $e$ indicates the patch that contains the point $\pt{x}_{k}$.
The matrix $[\hat{\mathbf{ A}}]$ and the vector $\{\bar{\mathbf{ c}}\} $ have to be modified for these grid points.

For the computation of the initial stresses we also require the strains. 
The strains are computed by:
\begin{equation}
\label{Bhat}
\{\myVecGreek{\epsilon}\}= [\hat{\mathbf{ B}}] \{\myVec{u}\}
\end{equation}
where $[\hat{\mathbf{ B}}]$ is presented later.

The values of initial stress at grid points are computed by:
\begin{equation}
\label{ }
\{\myVecGreek{\sigma}_{0}\}= [\mathbf{ D}^{\prime}] [\hat{\mathbf{ B}}]\{\myVec{u}\}
\end{equation}
where $[\mathbf{ D}^{\prime}] $ is  matrix that contains sub-matrices $\mathbf{ D}^{\prime}$ on the diagonal.

\section{Solution procedure}
Eqs. (\ref{DisIE}), (\ref{Bhat}) form a linear system of equations. In fact, the initial stresses are a function of the elastic strains. Therefore, the system of equations may be solved either iteratively or by avoiding iterations altogether (referred to as a one step solution here).

\subsection{Iterative solution for elastic effects}
To solve iteratively we first solve
\begin{equation}
\label{DisIE1}
[\mathbf{ L}] \{\mathbf{ x}\}_{0} = \{\mathbf{ r}\} 
\end{equation}
and then compute increments of the solution vector $\{ \triangle \mathbf{ x}\}_{i}$ due to the effect of the initial stress:
\begin{equation}
\label{DisIE2}
[\mathbf{ L}] \{\triangle\mathbf{ x}\}_{i} = [\mathbf{ B}_{0}] \{ \myVecGreek{\sigma}_{0}\} 
\end{equation}
where the subscript $i$ is the iteration number.
The final values are obtained by summing all the increments after convergence:
\begin{equation}
\label{Sumiter}
\{ \mathbf{ x}\}_{i} = \{ \mathbf{ x}\}_{0} + \{ \triangle \mathbf{ x}\}_{1} + \{ \triangle \mathbf{ x}\}_{2} \cdots
\end{equation}

\subsection{One step solution}
To avoid having to iterate for the difference in elastic properties we present an option of a solution that already includes this effect.
Eq. (\ref{Bhat}) can be written in the following form:
\begin{equation}
\label{Strain_rev1}
\{ \myVecGreek{\epsilon}\}= [\hat{\mathbf{ C}}]  \{ \mathbf{ x} \} + \{\bar{\bar{\mathbf{ c}}}\} + [\hat{\mathbf{ C}}_{0}]  [ \mathbf{ D}^{\prime}] \{ \myVecGreek{\epsilon} \}
\end{equation}

We define
\begin{equation}
\label{}
[\hat{\mathbf{ C}}] = [\hat{\mathbf{ B}}] [\hat{\mathbf{ A}}] \hspace{10mm}
[\hat{\mathbf{ C}}_{0}] = [\hat{\mathbf{ B}}] [\bar{\mathbf{B}}_{0}] \hspace{10mm}
\{\bar{\bar{\mathbf{ c}}}\} = [\hat{\mathbf{ B}}] \{\bar{\mathbf{ c}}\} 
\end{equation}
Eq. (\ref{Strain_rev1}) along with Eq. (\ref{DisIE}) form the following linear system of equations:
\begin{equation}
\label{Final_sys}
\begin{pmatrix}
  [\mathbf{ L}]        & - [\mathbf{ B}_{0}] [\mathbf{D}^{\prime}]  \\ \\
 - [\hat{\mathbf{ C}}] & [{\mathbf{I}}] - [\hat{\mathbf{ C}}_{0}] [\mathbf{D}^{\prime}]
\end{pmatrix}
\begin{pmatrix}
\{ \mathbf{ x}\} \\ \\
\{ \myVecGreek{\epsilon}\}
\end{pmatrix}
=
\begin{pmatrix}
\{\mathbf{ r}\} \\ \\
\{\bar{\bar{\mathbf{ c}}}\}
\end{pmatrix}
\end{equation}
that can be solved in terms of boundary unknowns and internal strains.

To avoid the construction of  Eq. (\ref{Final_sys}), it is possible to obtain a system of equations in terms of the boundary unknowns only:
\begin{equation}
\label{onestep}
[\mathbf{ L}]^{\prime} \{ \mathbf{ x}\} = \{\mathbf{ r}\}^{\prime} 
\end{equation}
where $[\mathbf{ L}]^{\prime}$ and $\{\mathbf{ r}\}^{\prime} $ are modified left and right hand sides that will be shown.

The strain vector can be computed from Eq. (\ref{Strain_rev1}):
\begin{equation}
\label{Strain_rev2}
\{ \myVecGreek{\epsilon}\}= ([\mathbf{I}] - 
[\hat{\mathbf{ C}}_0] [\mathbf{D}^{\prime}])^{-1} ([\hat{\mathbf{C}}] \{ \mathbf{ x}\} + \{\bar{\bar{\mathbf{ c}}}\}) = [\mathbf{A}] \{\mathbf{x}\} + \{\mathbf{b}\}
\end{equation}

\noindent where:

\begin{equation}
[\mathbf{A}] = ( [\mathbf{I}] - [\hat{\mathbf{C}}_0] [\mathbf{D}^{\prime}])^{-1} [\hat{\mathbf{C}}] \hspace{5mm} \{\mathbf{b}\} = ( [\mathbf{I}] - [\hat{\mathbf{C}}_0] 
[\mathbf{D}^{\prime}])^{-1} \{\bar{\bar{\mathbf{ c}}}\}
\end{equation}

The strain vector in Eq. (\ref{Strain_rev2}) can be replaced in Eq. (\ref{DisIE}) in order to obtain:
\begin{equation}
[\mathbf{L}] \{\mathbf{x}\} = \{\mathbf{r}\} + [\mathbf{B}_0] 
[\mathbf{D}^{\prime}] ([\mathbf{A}] \{\mathbf{x}\} + \{\mathbf{b}\})
\end{equation}
and, hence, the following system of equations to be solved in terms of $\{\mathbf{x}\}$:

\begin{equation}
\label{Solve}
( [\mathbf{L}] - [\mathbf{B}_0] [\mathbf{D}^{\prime}] [\mathbf{A}] ) 
\{\mathbf{x}\} = \{\mathbf{r}\} + [\mathbf{B}_0] 
[\mathbf{D}^{\prime}] \{\mathbf{b}\}
\end{equation}
The matrices in Eq. (\ref{onestep}) are defined by:
\begin{eqnarray}
[\mathbf{L}]^{\prime} & = &  [\mathbf{L}] - [\mathbf{B}_0] [\mathbf{D}^{\prime}] [\mathbf{A}] \\
\{\mathbf{ r}\}^{\prime} & = &\{\mathbf{r}\} + [\mathbf{B}_0] 
[\mathbf{D}^{\prime}] \{\mathbf{b}\} 
\end{eqnarray}

\section{Computation of $\mathbf{ B}_0$-matrix}
Since we are dealing here with thin inclusions the following assumptions are made:
\begin{itemize}
  \item The variation of displacements across the thickness is either linear or constant.
  \item The shear stresses across the thickness are neglected.
\end{itemize}
For thin inclusions it is convenient and efficient to assume that the initial stresses $ \myVecGreek{\sigma}_{0}$ are defined in local directions ($x^{\prime},y^{\prime},z^{\prime}$), especially if they are curved.
The volume integral is changed to:
\begin{equation}
\mathbf{ B}^{\prime m}_{0n}  =  \int_{\domain_m}  \fund{E}^{\prime} (\sourcept_{n},\fieldpt)   \myVecGreek{\sigma}^{\prime}_{0} (\fieldpt) d \domain_m (\fieldpt) 
\end{equation}
where $\domain_m$ defines the volume of inclusion $m$ and  $\myVecGreek{\sigma}^{\prime}_{0} (\fieldpt)$ specify the local initial stresses at point $\fieldpt$ inside the inclusion.

To compute $\fund{E}^{\prime}$ we consider that the index $i$ relates to the direction of the source and the indexes $j,k$ to the strain component and that only the strain is required to be computed in the local direction.  
The fundamental solution $\fund{E}$ in global directions is given for example for $i=1$:
\begin{equation}
\label{ }
E_{1jk}= \frac{-C}{r^{2}}\left[C_{3}(r_{,k} \ \delta_{1j} + r_{,j} \ \delta_{1k}) - r_{,1} \ \delta_{jk} + C_{4} \ r_{,1} \ r_{,j}  \ r_{,k}\right]
\end{equation}
To convert the strain into local directions we make the following transformation:
\begin{equation}
\label{ }
\mathbf{ E}^{\prime}_{1}= \mathbf{ R}^{T} \mathbf{ E}_{1} \mathbf{ R}
\end{equation}
where $\mathbf{ E}_{1}$ is a tensor containing components $E_{1jk}$ and $\mathbf{ R}$ is the transformation tensor given by:
\begin{equation}
\label{ }
\mathbf{ R}= \left(\begin{array}{ccc}\mathrm{v}_{1_{x}} & \mathrm{v}_{2_{x}} & \mathrm{v}_{3_{x}} \\\mathrm{v}_{1_{y}} & \mathrm{v}_{2_{y}} & \mathrm{v}_{3_{y}} \\\mathrm{v}_{1_{z}} & \mathrm{v}_{2_{z}} & \mathrm{v}_{3_{z}}\end{array}\right)
\end{equation}
where $\mathbf{ v}_{1}, \mathbf{ v}_{2}, \mathbf{ v}_{3}$ are unit vectors in the local coordinate directions (see section \ref{Geom}).
We do this transformation for all values of $i$ and then introduce the Voight notation.

\section{Numerical volume integration}
We use Gauss Quadrature for evaluating the volume integral. However, we have to consider that the integrand tends to infinity with $O(r^2)$ as point $\sourcept$ is approached.
We subdivide the inclusion region into integration regions, which are determined depending on the location of points $\sourcept$ and the aspect ratio of the sub-region.
We then have to consider two cases: One where the point $\sourcept$ is on the edge of an integration region (singular integration) and one where it is not (regular integration).

With respect to singular integration we explore two options:
\begin{itemize}
  \item Option 1: Numerical integration
  \item Option 2: Combined numerical and analytical integration
\end{itemize}

\subsection{Regular integration}
For integration region $n_{s}$ the transformation from the coordinates used for Gauss integration $\bar{\myVecGreek{\xi}}=(\bar{\xi},\bar{\eta},\bar{\zeta})^{\mathrm{T}}=[-1,1]^3$ to the ($\pt{s}=(s,t,r)^{\mathrm{T}}=[0,1]^3$) coordinates of the inclusion is  given by
\begin{eqnarray}
\label{Gauss}
s & = \frac{\Delta s_{n}}{2} (1+\bar{\xi}) + s_{1n} \\
\nonumber
t & = \frac{\Delta t_{n}}{2} (1+\bar{\eta}) + t_{1n} \\
\nonumber
r & = \frac{\Delta r_{n}}{2} (1+\bar{\eta}) + r_{1n} 
\end{eqnarray}
where $\Delta s_{n}\times \Delta t_{n},\times \Delta r_{n} $ denotes the size of the integration region and $s_{1n},t_{1n},r_{1n}$ are the coordinates of the lower left edge.
The Jacobian of this transformation is $J_{\xi}^{n_{s}}=\frac{1}{8}\ \Delta s_{n} \  \Delta t_{n} \  \Delta r_{n}$.

Leaving out the superscript $m$ that denotes the inclusion number, we can write:
\begin{equation}
  \label{Integ3D}
   \mathbf{ B}^{\prime}_{0n} =  \sum_{n_{s}=1}^{N_{s}}\int_{-1}^{1} \int_{-1}^{1} \int_{-1}^{1}  \fund{E}^{\prime} \left( \sourcept_{n},\bar{\pt{x}}(\bar{\xi},\bar{\eta},\bar{\zeta}) \right)
\myVecGreek{\sigma}^{\prime}_{0} \left( \bar{\pt{x}} (\bar{\xi},\bar{\eta},\bar{\zeta}) \right) J(\pt{s}) \ J_{\xi}^{n_{s}} \ d \bar{\xi} d \bar{\eta} d \bar{\zeta} 
\end{equation}
where $ J(\mathbf{ s})$ is the Jacobian of the mapping between $\pt{s}$ and $\pt{x}$ coordinate systems.

Applying Gauss integration we have:
\begin{eqnarray}
  \label{Gauss3D}
   \mathbf{B}^{\prime}_{0n}  \approx  \sum_{n_{s}=1}^{N_{s}} \sum_{g_{s}=1}^{G_{s}} \sum_{g_{t}=1}^{G_{t}}   \sum_{g_{r}=1}^{G_{r}} \fund{E}\left( \sourcept_{n},\bar{\pt{x}}(\bar{\xi}_{g_{s}},\bar{\eta}_{g_{t}},\bar{\zeta}_{g_{r}}) \right)\myVecGreek{\sigma}^{\prime}_{0} \left( \bar{\pt{x}}(\bar{\xi}_{g_{s}},\bar{\eta}_{g_{t}},\bar{\zeta}_{g_{r}}) \right) \\
 \nonumber
    J(\pt{s})  \ J_{\xi}^{n_{s}} \ W_{g_{s}} \ W_{g_{t}} \ W_{g_{r}} 
\end{eqnarray}
where $N_{s}$ is the number of subregions and $G_{s},G_{t},G_{r}$  are the number of Gauss points (which depends on the proximity of $\sourcept_n$ to the region) and $\bar{\xi}_{g_{s}},\bar{\eta}_{g_{t}},\bar{\zeta}_{g_{t}}$ the Gauss point coordinates in $s, t,r$  directions, respectively. $W_{g_{s}},W_{g_{t}},W_{g_{r}}$ are Gauss weights. 

\subsection{Singular integration}
When $\sourcept$ is on one of the edges of an integration region the integrand becomes singular.
To perform the integration, we define a zone of exclusion $\domain_{e}$ at $\sourcept$ and split the integral into two parts:
\begin{equation}
\mathbf{ B}^{\prime}_{0n}  =  \int_{\domain_{m} - \domain_{e}}  \fund{E}^{\prime} (\sourcept_{n},\fieldpt)   \myVecGreek{\sigma}^{\prime}_{0} (\fieldpt) d \domain  + \int_{\domain_{e}}  \fund{E}^{\prime} (\sourcept_{n},\fieldpt)   \myVecGreek{\sigma}^{\prime}_{0} (\fieldpt) d \domain_{e} 
\end{equation}
The first part is nearly singular and can be integrated using Gauss Quadrature. The second part is singular and can be integrated using Gauss Quadrature with special mapping procedures or integrated analytically.

\subsubsection{Integration of nearly singular part}
\begin{figure}
\begin{center}
\begin{overpic}[scale=0.4]{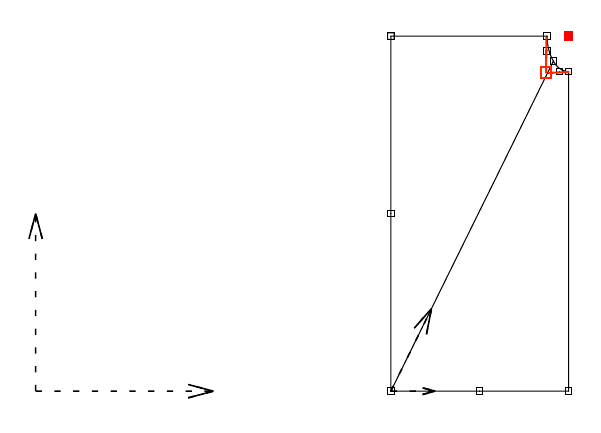} 
\put(70,22){$\hat{\vv}$}
\put(70,10){$\hat{\uu}$}
\put(35,10){$s$}
\put(10,35){$t$}
\put(10,10){$r$}
\put(95,35){$L_{t}$}
\put(75,70){$L_{s}$}
\put(95,70){$\sourcept$}
\put(96,62){$R_{2}$}
\put(91,68){$R_{1}$}
\put(60,2){$s_1$}
\put(90,2){$s_2$}
\put(60,10){$t_1$}
\put(60,67){$t_2$}
\end{overpic}
\caption{View of the subdivision of integration region containing $\sourcept$ into quadratic-linear subregions, in the local $s,t,r$ coordinate system. Control points are shown as hollow squares. For the linear-linear basis functions we only use the control points at the corners and a changed control point marked in red. }
\label{Singular1L}
\end{center}
\end{figure}
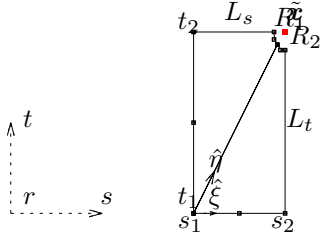

We divide the region of integration into two subregions which are tapered towards the singular point, excluding a zone containing $\sourcept$ (\myfigref{Singular1L}).
Two mappings are required. One from the $\bar{\myVecGreek{\xi}}=(\bar{\xi},\bar{\eta},\bar{\zeta})^{\mathrm{T}}=[-1,1]^3$coordinate system (where the Gauss points are defined) to coordinates $\hat{\myVecGreek{\xi}}=(\hat{\xi},\hat{\eta},\hat{\zeta})^{\mathrm{T}}=[0,1]^3$:
\begin{eqnarray}
\label{MapG}
\hat{\uu}= 0.5(1+\bar{\uu}) \\
\hat{\vv}= 0.5(1+\bar{\vv}) \\
\hat{\zeta}= 0.5(1+\bar{\zeta}) 
\end{eqnarray}
The Jacobian of this transformation is $J_{\bar{\xi}}=(0.5)^3$.

The second mapping is from the $\hat{\myVecGreek{\xi}}$ to the local ($\pt{s}$) coordinate system of the inclusion :
\begin{eqnarray}
\label{ }
\left\{\begin{array}{c}s \\ t \end{array}\right\} &=& \sum^{I}_{i=1} \sum^{J}_{j=1} B_{ij}(\hat{\uu},\hat{\vv}) \left\{\begin{array}{c} s_{ij} \\ t_{ij} \end{array}\right\} \\
r&=& \hat{\zeta}
\end{eqnarray}
where $B_{ij}$ are NURBS basis functions, $I,J$ are the number of control points in each direction and $s_{ij}, t_{ij} $ are control point coordinates. 
Two choices for $B_{ij}$ are possible: linear-linear and quadratic-linear. The first choice results in a prismatic the second in a cylindrical region of exclusion.

Table \ref{tab:Cm} shows how the control points are defined for a quadric-linear NURBS and for the example in  \myfigref{Singular1L}.  We define $R_1= R/L_s$ and $R_2= R/L_t$ where $R$ is the radius of the exclusion zone and $L_{s}, L_{t}$ the size of the integration region in $s,t$ directions as defined in \myfigref{Singular1L}. $s_1,s_2,t_1,t_2$ are the edge coordinates of the integration region, where we drop the subscript $n$ designating the region number. $s_p,t_p$ are the local coordinates of the singular point.
\begin{mytable}
  {H}               
  {Control point coordinates and weights for the singular point on right upper corner}  
  {tab:Cm}  
  {cccccccc}         
  \mytableheader{ Sub-region &  ij= & 11 & 21 & 31 & 12 & 22 & 32 }  
1&$ s_{ij}$=& $s_{1}$&  $0.5(s_{1}+ s_{p})$ & $s_{p}$ & $s_{p} -0.707R_1$ & $s_{p}-0.414R_1$ & $ s_{p}$\\ \\
&$ t_{ij}$=& $t_{1}$&  $t_{1}$ & $t_{1}$ & $t_{p} - 0.707R_2$ & $t_{p}-R_2$ & $ t_p - R_2$\\ \\
&$ w_{ij}$=& 1&  1 & 1 & 1 & 0.924 & $ 1$\\ \\
2&$ s_{ij}$=& $s_{1}$&  $s_{1}$ & $s_{1}$ & $s_{p} - R_1$ & $s_{p} - R_1$ & $ s_{p} - 0.707R_1$\\ \\
&$ t_{ij}$=& $t_{p}$&  $0.5(t_{1}+ t_{p})$ & $t_{1}$ & $t_{p} $ & $t_{p}- 0.414R_2$ & $ t_{p} - 0.707R_2$\\ \\
&$ w_{ij}$=& 1&  1 & 1 & 1 & 0.924 & $ 1$\\
\end{mytable}%

\textbf{The Jacobian of this mapping is $J_{\xi}^{n_{s}}$ and tends  with \textit{O}(r) to zero as $\sourcept$ is approached. 
This means that the singularity of the integrand is reduced by one order and the integrand becomes weakly singular.}
However, since $\sourcept$ is excluded the integrand remains nearly singular and can be integrated using Gauss Quadrature with careful selection of the number of Gauss points.
The final mapping is from the $\pt{s}$ to the $\pt{x}$ coordinate system with Jacobian $J(\mathbf{ s})$.

The numerical integration of a nearly singular subregion is given by
\begin{eqnarray}
  \label{}
  \mathbf{B}_{0n}  \approx  \sum_{n_{t}=1}^{N_{t}} \sum_{g_{s}=1}^{G_{s}} \sum_{g_{t}=1}^{G_{t}}   \sum_{g_{r}=1}^{G_{r}} \fund{E}\left( \sourcept_{n},\bar{\pt{x}}(\bar{\xi}_{g_{s}},\bar{\eta}_{g_{t}},\bar{\zeta}_{g_{r}}) \right)\myVecGreek{\sigma}^{\prime}_{0} \left( \bar{\pt{x}}(\bar{\xi}_{g_{s}},\bar{\eta}_{g_{t}},\bar{\zeta}_{g_{r}}) \right) \\
 \nonumber
    J(\pt{s})  \ J_{\xi}^{n_{s}} \ J_{\bar{\xi}} \ W_{g_{s}} \ W_{g_{t}} 
\end{eqnarray}
where $N_{t}$ is the number of triangular subregions.

\newpage

\subsubsection{Numerical integration of singular part}
For this we divide the prismatic exclusion region into  tetrahedral sub-regions. Within a tetrahedral sub-region the Jacobian tends to zero with $O(r^2)$ as the singular point is approached, cancelling out the singularity.
\begin{figure}
\begin{center}
\begin{overpic}[scale=0.5]{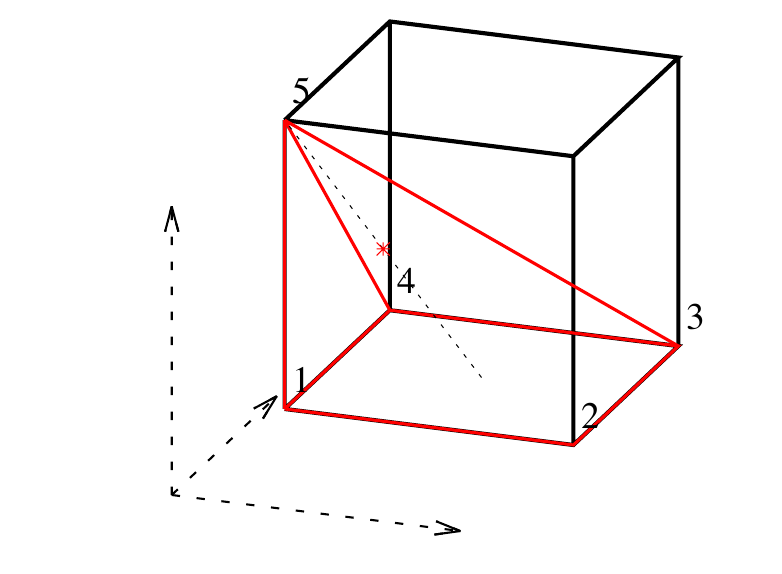}
 \put(60,5){$s$}
  \put(24,45){$r$}
    \put(30,25){$t$}
\end{overpic}
\begin{overpic}[scale=0.5]{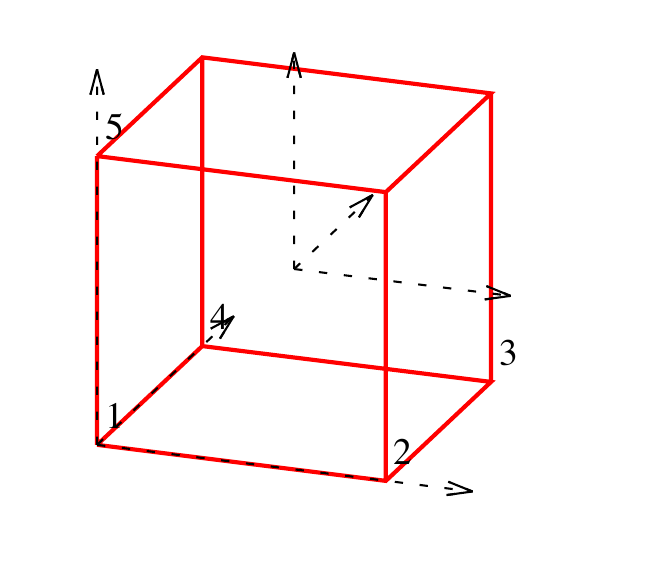}
\put(75,10){$\sigma$}
\put(10,75){$\rho$}
\put(37,37){$\tau$}
\put(80,40){$\bar{\xi}$}
\put(55,60){$\bar{\eta}$}
\put(45,75){$\bar{\zeta}$}
\end{overpic}
\caption{Singular volume integration, showing a tetrahedral subregion and the mapping from the $\pt{s}$ to the 
$ \sigma,\tau,\rho$  coordinate system.}
\label{Tetra}
\end{center}
\end{figure}
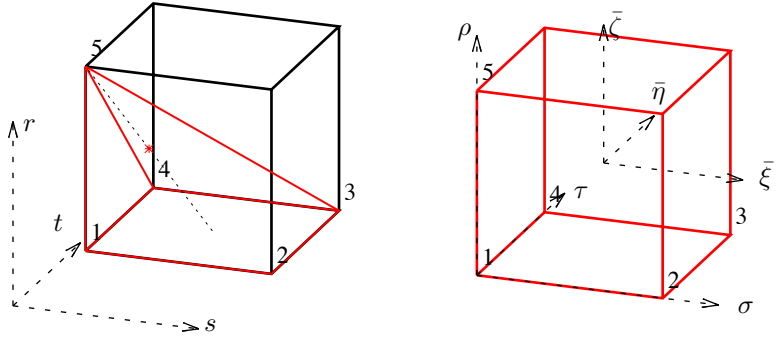

The transformation from the local $\myVecGreek{\bar{\xi}}$ coordinate system, in which the Gauss coordinates are defined, to global coordinates involves the following transformation steps:
\begin{enumerate}
  \item from $\myVecGreek{\bar{\xi}}$ to a local tetrahedral system $(\sigma,\tau,\rho)^{\mathrm{T}}=[0,1]^3$ 
  \item from ($\sigma,\tau,\rho$) to $\pt{s}$
  \item from $\pt{s}$ to $\pt{x}$
\end{enumerate}
 
The first transformation is:
\begin{eqnarray}
\label{}
\sigma= 0.5(1+\bar{\uu}) \\
\tau= 0.5(1+\bar{\vv}) \\
\rho= 0.5(1+\bar{\zeta}) 
\end{eqnarray}
The Jacobian of this transformation is $J_{\bar{\xi}}=(0.5)^3$.

Referring to \myfigref{Tetra} we assume that the singular point is an edge point of the integration region.
For this case the transformation is as follows:
First we determine the local coordinates $\pt{ s}_{1}$ to $\pt{ s}_{5}$ of the edge points of the tetrahedron, with 5 being the singularity point.
Next we define a linear plane NURBS surface with points 1 to 4 and map the coordinates of the point ($\sigma, \tau$) onto this surface:
\begin{equation}
\pt{ s}_{0}(\sigma,\tau)= \sum_{i=1}^{4}R_{i}(\sigma,\tau) \ \pt{ s}_{i}
\end{equation}
where $R_{i}(\sigma,\tau)$ are linear basis functions.
The final map is obtained by a linear interpolation in the $\rho$-direction:
\begin{equation}
\pt{ s}(\sigma,\tau,\rho)= (1-\rho) \ \mathbf{ s}_{0}(\sigma,\tau) + \rho \  \pt{ s}_{5}
\end{equation}
The Jacobi matrix of this transformation is given by:
\begin{equation}
\mathbf{ J}= \left(\begin{array}{c}(1-\rho)\frac{\partial \pt{ s}_{0}}{\partial\sigma} \\ \\ (1-\rho)\frac{\partial \pt{ s}_{0}}{\partial\tau} \\ \\ \pt{ s}_{5} - \pt{ s}_{0}\end{array}\right)
\end{equation}
The Jacobian of this transformation $J_{t}= |\mathbf{ J}|$  tends to zero with $O(r^2)$ as the singular point ($\rho=1$) is approached.

The numerical integration is:
\begin{eqnarray}
  \label{}
  \mathbf{B}_{0n}  \approx  \sum_{n_{t}=1}^{N_{t}} \sum_{g_{s}=1}^{G_{s}} \sum_{g_{t}=1}^{G_{t}}   \sum_{g_{r}=1}^{G_{r}} \fund{E}\left( \sourcept_{n},\bar{\pt{x}}(\bar{\xi}_{g_{s}},\bar{\eta}_{g_{t}},\bar{\zeta}_{g_{r}}) \right)\myVecGreek{\sigma}^{\prime}_{0} \left( \bar{\pt{x}}(\bar{\xi}_{g_{s}},\bar{\eta}_{g_{t}},\bar{\zeta}_{g_{r}}) \right) \\
 \nonumber
    J(\pt{s})  \ J_{t}^{n_{t}} \ J_{\bar{\xi}} \ W_{g_{s}} \ W_{g_{t}} 
\end{eqnarray}
where $N_{t}$ is the number of tetrahedral subregions.

\subsubsection{Analytical integration of singular part}
As the thickness of the inclusion decreases the error in the numerical integration of the region of exclusion increases, as the aspect ratio of the tetrahedral subregions  becomes extreme. In this case we can switch to an analytical integration of a cylindrical exclusion region.
For this we assume that  $\myVecGreek{\sigma}^{\prime}_{0}$ is constant inside the exclusion region, i.e. we approximate the integral by:
\begin{equation}
\int_{\domain_{e}}  \fund{E}^{\prime} (\sourcept_{n},\fieldpt)   \myVecGreek{\sigma}^{\prime}_{0} (\fieldpt) d \domain_{e} \approx \left( \ \int_{\domain_{e}}  \fund{E}^{\prime} (\sourcept_{n},\fieldpt)   d \domain_{e} \right) \myVecGreek{\sigma}^{\prime}_{0} (\sourcept_n)
\end{equation}

The integral in the parentheses can now be evaluated analytically.
\begin{figure}
\begin{center}
\begin{overpic}[scale=0.25]{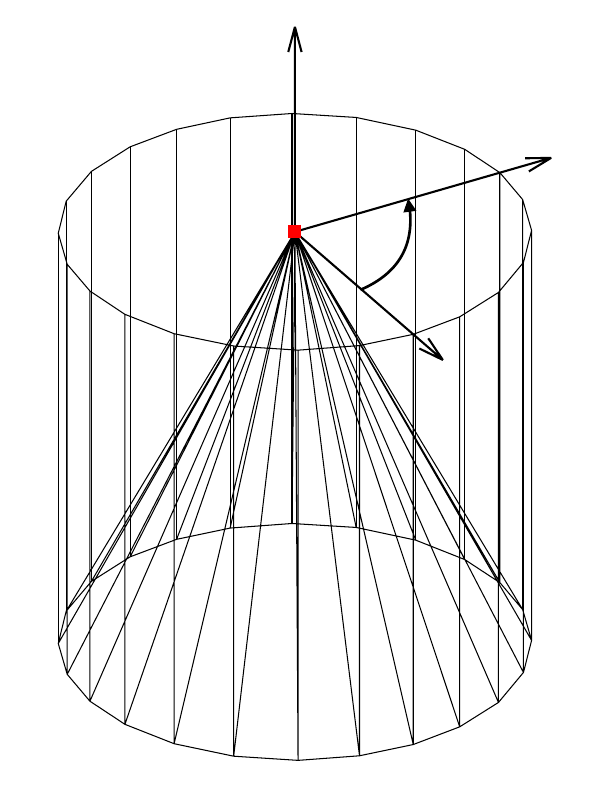}
 \put(40,15){$\phi$}
 \put(7,10){$x^{\prime}$}
 \put(80,15){$y^{\prime}$}
 \put(40,95){$z^{\prime}$}
 \put(48,50){$r$}
  \put(40,82){$\theta$}
 \put(10,50){$H$}
 \put(50,85){$R$}
\end{overpic}
\caption{Analytical integration regions for the case where the singular point is at the top.}
\label{Anal}
\end{center}
\end{figure}
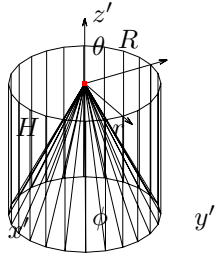
For this we define a cylindrical domain and split it into 2 regions shown in \myfigref{Anal}.
In cylindrical coordinates the integral becomes:

\begin{equation}
\triangle \fund{E}^{\prime} = \triangle \fund{E}^{\prime}_1+\triangle \fund{E}^{\prime}_2
\end{equation}

For the case where the singular point is at the top, the above terms are given by:

\begin{equation}
\label{I1I2}
\triangle \fund{E^{\prime}}_1  =  \int_{\phi_1}^{\phi_2} \int_{\theta=\pi-\tilde{\theta}}^{\pi} \int_{r=0}^{\frac{H}{\cos({\pi-\theta})}}  \frac{1}{r^{2}} \fund{\tilde{E}^{\prime}} \sin{\theta}  dr \ r^{2} d\theta d\phi 
\end{equation}
\begin{equation}
\triangle \fund{E^{\prime}}_2  =  \int_{\phi_1}^{\phi_2} \int_{\theta=\pi/2}^{\pi-\tilde{\theta}} \int_{r=0}^{\frac{R}{\sin({\pi-\theta})}} \frac{1}{r^{2}}  \fund{\tilde{E}^{\prime}} \sin{\theta}  dr \ r^{2} d\theta d\phi \nonumber
\end{equation}
where $\phi$ is measured counterclockwise from the $x^{\prime}$ axis.

For the case where the singular point is at the bottom, we have:

\begin{equation}
\label{I1I2}
\triangle \fund{E^{\prime}}_1  =  \int_{\phi_1}^{\phi_2} \int_{\theta=0}^{\tilde{\theta}} \int_{r=0}^{\frac{H}{\cos{\theta}}}  \frac{1}{r^{2}} \fund{\tilde{E}^{\prime}} \sin{\theta}  dr \ r^{2} d\theta d\phi 
\end{equation}
\begin{equation}
\triangle \fund{E^{\prime}}_2  =  \int_{\phi_1}^{\phi_2} \int_{\theta=\tilde{\theta}}^{\pi/2} \int_{r=0}^{\frac{R}{\sin{\theta}}} \frac{1}{r^{2}}  \fund{\tilde{E}^{\prime}} \sin{\theta}  dr \ r^{2} d\theta d\phi \nonumber
\end{equation}
with $\tilde{\theta}= \arctan (R/H)$. It can be seen that the $r^{2}$ terms cancel out which means that the integrand is no longer singular.

The analytical solutions $\triangle \fund{E}_{bottom}^{\prime}(i,j)$ for the range of $\phi= \phi_1 : \phi_2$ for the case where the singular point is at the top are given in the Appendix. The solutions for the singular point at the bottom are given by:

\begin{equation}
\triangle \fund{E}_{bottom}^{\prime}(i,j)  = - \triangle \fund{E}_{top}^{\prime}(i,j)
\mbox{ for   }  i=1,2 \mbox{  and  } j=5,6
\mbox{  and for    }  i=3 \mbox{  and  } j=1:4
\end{equation}

A transformation to the global system is necessary:
\begin{equation}
\label{ }
\mathbf{ B}_{0nj}= \mathbf{ T} \mathbf{ B}_{0nj}^{\prime}
\end{equation}

\newpage

\section{Computation of local strains inside the inclusion}
Local strains at grid points are computed from the values of displacements $\myVec{\primary}$ at grid points, using derivatives of interpolation functions $M_{k}$.
The interpolation functions can be either Lagrange polynomials or NURBS.
The advantage of using NURBS is that the limitations of Lagrange polynomials can be overcome (this will be discussed in more detail for the third test example). However, since NURBS use parameters instead of real values, care has to be taken.
In the following we first deal with the interpolation using Lagrange polynomials and then introduce NURBS.

\begin{figure}
\begin{center}
\begin{overpic}[scale=0.3]{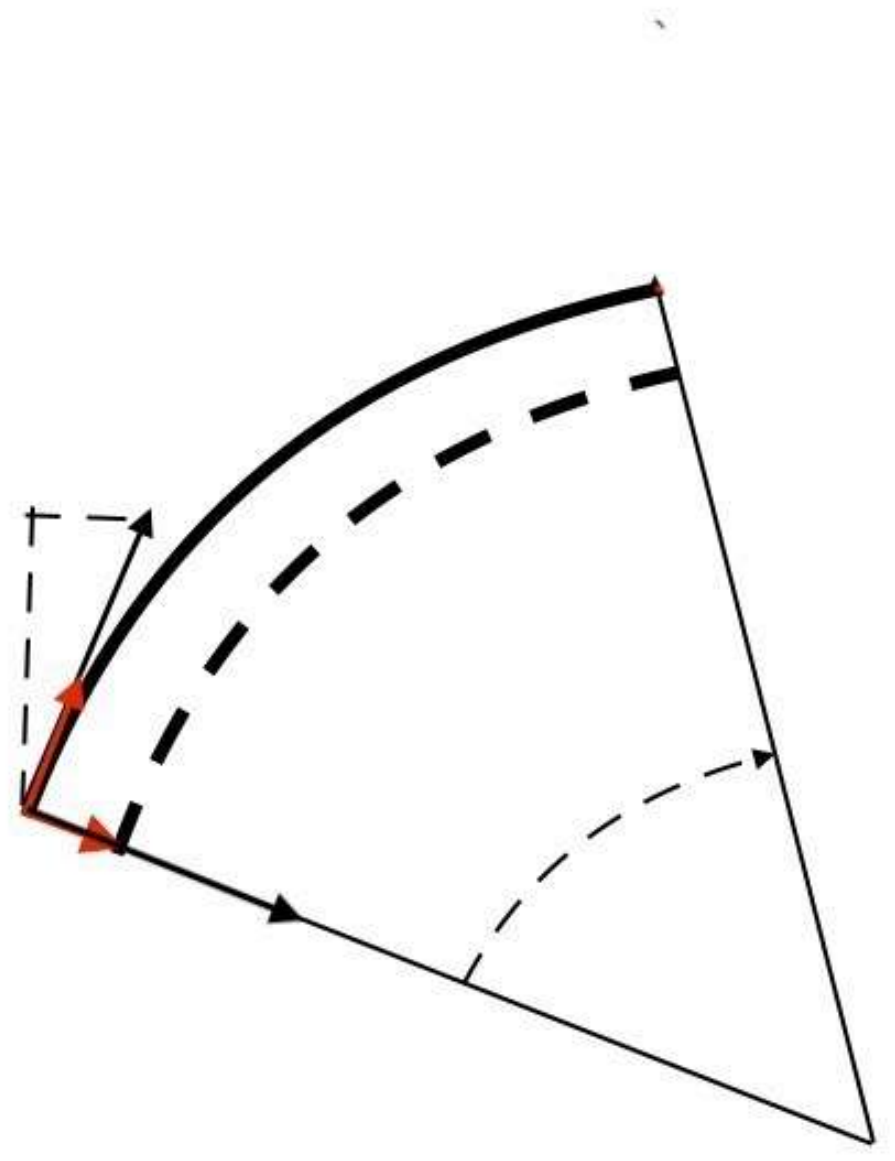}
 \put(12,27){$u_{z^{\prime}}$}
 \put(6,40){$u_{x^{\prime}}$}
 \put(22,60){$\mathbf{ V}_{s}$}
  \put(30,30){$\mathbf{ V}_{r}$}
\put(65,20){$d\phi$}
\put(30,65){$dx^{\prime}$}
\put(30,45){$dx^{\prime} - du^{\kappa}_{x^{\prime}}$}
\put(70,40){$R$}
\end{overpic}
\caption{Explanation of the computation of $\epsilon^{\kappa}_{x^{\prime} x^{\prime}}$ due to curvature.}
\label{Strain}
\end{center}
\end{figure}
We explain the strain computation first for the strain in the $x^{\prime}$ direction.
The local strain $\epsilon_{x^{\prime} x^{\prime}}$ has two components. One resulting from a change in displacement in $x'$-direction (which corresponds to the local $s$-direction):
\begin{equation}
\label{ }
\epsilon^s_{x^{\prime} x^{\prime}}=  \frac{\partial \primary_{x^{\prime}}}{\partial x^{\prime}} =\frac{\partial \primary_{x^{\prime}}}{\partial s}  \frac{\partial s}{\partial x^{\prime}}
\end{equation}
The other resulting from the curvature. Referring to \myfigref{Strain} the length of a small segment of an inclusion surface (in the $s,r$ plane) is given by:
\begin{equation}
\label{ }
dx^{\prime}= Rd \phi
\end{equation}
The changed length due to curvature and $u_{z^{\prime}}$ is 
\begin{equation}
\label{ }
dx^{\prime} - du^{\kappa}_{x^{\prime}}= (R - u_{z^{\prime}}) d \phi
\end{equation}
Therefore the change in length due to curvature is given by:
\begin{equation}
\label{ }
du^{\kappa}_{x^{\prime}}= -u_{z^{\prime}}  d\phi
\end{equation}
and the strain is:
\begin{equation}
\label{ }
\epsilon^{\kappa}_{x^{\prime} x^{\prime}}= \frac{du^{\kappa}_{x^{\prime}}}{dx^{\prime}} = -u_{z^{\prime}}  \frac{d \phi}{d x^{\prime}}
\end{equation}
Substitution of $\kappa_{x^{\prime}}= \frac{d \phi}{d x^{\prime}}$, where $\kappa_{x^{\prime}}$ is the curvature along $x^{\prime}$, we have
\begin{equation}
\label{ }
\epsilon^{\kappa}_{x^{\prime} x^{\prime}}=   -u_{z^{\prime}}  \kappa_{x^{\prime}}
\end{equation}
The curvature is computed by
\begin{equation}
\label{ }
\kappa_{x^{\prime}}= \left| \frac{\partial \pt{x}}{\partial s} \times  \frac{\partial^2 \pt{x}}{\partial s^2} \right| \frac{1}{J_{s}^3}
\end{equation}
where $J_{s}= |\mathbf{ V}_{s} |$.

The local strains in $x^{\prime}, y^{\prime} z^{\prime}$ directions are given by:
\begin{eqnarray}
\label{eq9:strain}
\epsilon_{x^{\prime} x^{\prime}}& = & \frac{\partial \primary_{x^{\prime}}}{\partial x^{\prime}} -  u_{z^{\prime}} \kappa_{x^{\prime}}\\
\epsilon_{y^{\prime} y^{\prime}} & = & \frac{\partial \primary_{y^{\prime}}}{\partial y^{\prime}} - u_{z^{\prime}} \kappa_{y^{\prime}} \\
\epsilon_{z^{\prime} z^{\prime}} & = & \frac{\partial \primary_{z^{\prime}}}{\partial z^{\prime}}  \\
\epsilon_{x^{\prime}y^{\prime}} &= & \frac{\partial \primary_{y^{\prime}}}{\partial x^{\prime}} + \frac{\partial \primary_{x^{\prime}}}{\partial y^{\prime}} 
\end{eqnarray}
where $ \kappa_{y^{\prime}}$ is the curvature in the $y^{\prime}$ direction.

\

\remark{It should be noted that the formula for local strain computations published in \cite{Gao2011} and \cite{BeerMarussig} does not include the term due to curvature and is therefore only valid for plane surfaces. To our best knowledge this is the first time the complete formula has been published.}

\

The local derivatives of the displacements are computed by
\begin{eqnarray}
\frac{\partial \primary_{x^{\prime}}}{\partial x^{\prime}} &=& \frac{\partial \primary_{x^{\prime}}}{\partial s} \frac{\partial s}{\partial x^{\prime}} \\
\frac{\partial \primary_{y^{\prime}}}{\partial y^{\prime}} &=& \frac{\partial \primary_{y^{\prime}}}{\partial s}\frac{\partial s}{\partial y^{\prime}} + \frac{\partial \primary_{y^{\prime}}}{\partial t} \frac{\partial t}{\partial y^{\prime}} \\
\frac{\partial \primary_{x^{\prime}}}{\partial y^{\prime}} &=& \frac{\partial \primary_{x^{\prime}}}{\partial s} \frac{\partial s}{\partial y^{\prime}} + \frac{\partial \primary_{x^{\prime}}}{\partial t} \frac{\partial t}{\partial y^{\prime}}  \\
\frac{\partial \primary_{y^{\prime}}}{\partial x^{\prime}}&=& \frac{\partial \primary_{y^{\prime}}}{\partial s} \frac{\partial s}{\partial x^{\prime}} \\
\frac{\partial \primary_{z^{\prime}}}{\partial z^{\prime}}&=&  \frac{\partial \primary_{z^{\prime}}}{\partial r} \frac{\partial r}{\partial z^{\prime}}
\end{eqnarray}
The displacements in local directions are given by:
\begin{eqnarray}
\label{ }
\primary_{x^{\prime}}= \myVec{\primary} \cdot \mathbf{v}_1 \\
\primary_{y^{\prime}}= \myVec{\primary} \cdot \mathbf{v}_2 \\
\primary_{z^{\prime}}= \myVec{\primary} \cdot \mathbf{v}_3 
\end{eqnarray}
where $ \mathbf{v}_1, \mathbf{v}_2, \mathbf{v}_3$ unit vectors in $x^{\prime},y^{\prime},z^{\prime}$ directions.

\newpage

\subsection{Determination of geometric derivatives}
\begin{figure}[H]
\begin{center}
\begin{overpic}[scale=0.35]{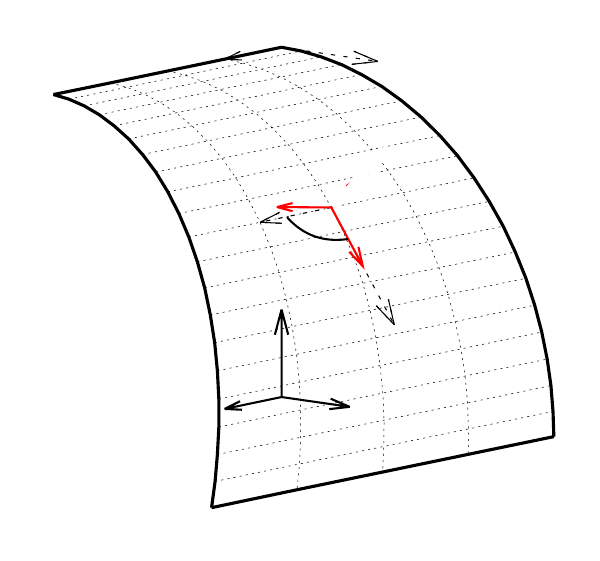}
 \put(62,48){$x^{\prime}$}
 \put(43,60){$y^{\prime}$}
 \put(43,52){$\mathbf{ v}_{t}$}
 \put(65,35){$\mathbf{ v}_{s}$}
 \put(53,55){$\theta$}
 \put(65,85){$s$}
 \put(40,85){$t$}
 \put(32,25){$x$}
 \put(60,25){$y$}
 \put(40,40){$z$}
\end{overpic}
\caption{Explanation for the determination of geometric derivatives.}
\label{Regularisation}
\end{center}
\end{figure}
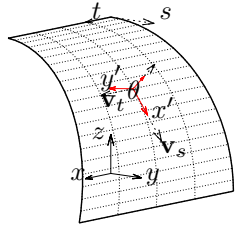

The geometric derivatives are given by:
\begin{eqnarray}
\frac{\partial s}{\partial x^{\prime}} & = & \frac{1}{J_{s}} \\
\frac{\partial s}{\partial y^{\prime}} & = & -\frac{\cos \theta}{J_{s} \ \sin \theta} \\
\frac{\partial t}{\partial x^{\prime}} & = & 0 \\
\frac{\partial t}{\partial y^{\prime}} & = & \frac{1}{J_{t}\ \sin \theta} \\
\frac{\partial r}{\partial z^{\prime}} & = & \frac{1}{J_{r}} 
\end{eqnarray}

We define
\begin{eqnarray}
J_{s}= |\mathbf{ V}_{s} |\\
J_{t}= |\mathbf{ V}_{t} | \\
J_{r}= |\mathbf{ V}_{r} |
\end{eqnarray}
where $\mathbf{ V}_{s} , \mathbf{ V}_{t} , \mathbf{ V}_{r} $ are vectors in $s,t,r$ directions.
and
\begin{eqnarray}
\cos \theta & = & \mathbf{ v}_{s} \cdot \mathbf{ v}_{t}\\
\sin \theta & = & \mathbf{ v}_{t} \cdot \mathbf{ v}_{y^{\prime}}
\end{eqnarray}

\subsection{Using Lagrange polynomials for interpolation}

Introducing the interpolation functions we obtain:
\begin{align}
\frac{\partial \primary_{x^{\prime}}}{\partial s}=\sum_{k=1}^{K}    \frac{\partial M_{k}}{\partial s}\  \primary_{x^{\prime}k}^{e} =\sum_{k=1}^{K}    \frac{\partial M_{k}}{\partial s}\  \left(\myVec{\primary} \cdot \mathbf{v}_{1} \right)_{k}^{e}\\
\frac{\partial \primary_{y^{\prime}}}{\partial t}=\sum_{k=1}^{K}    \frac{\partial M_{k}}{\partial t}\  \primary_{y^{\prime}k}^{e} =\sum_{k=1}^{K}    \frac{\partial M_{k}}{\partial t}\  \left(\myVec{\primary} \cdot \mathbf{v}_{2} \right)_{k}^{e}\\
\frac{\partial \primary_{z^{\prime}}}{\partial r}=\sum_{k=1}^{K}    \frac{\partial M_{k}}{\partial r}\  \primary_{z^{\prime}k}^{e} =\sum_{k=1}^{K}    \frac{\partial M_{k}}{\partial r}\  \left(\myVec{\primary} \cdot \mathbf{v}_{3} \right)_{k}^{e}
\end{align}
where $M_{k}(s,t,r)$ are Lagrange polynomials.
 
 \paragraph{Orthogonal local axes} 

If $ \mathbf{ v}_{t}$ is orthogonal to $ \mathbf{ v}_{s}$  then $\cos \theta = 0$ and we have:
\begin{eqnarray}
\frac{\partial s}{\partial x^{\prime}} & = & \frac{1}{J_{s}} \\
\frac{\partial s}{\partial y^{\prime}} & = &0 \\
\frac{\partial t}{\partial y^{\prime}} & = & \frac{1}{J_{t}}
\end{eqnarray}

In matrix algebra the strains at an internal point $i$ are :
\begin{equation}
\label{ }
 \myVecGreek{\epsilon}_{i}^{\prime}= \sum \mathbf{\hat{ B}}^{\prime}_{ij}\mathbf{ u}_{j}
\end{equation}
with
\begin{equation}
\label{ }
\mathbf{ \hat{B}}^{\prime}_{ij}=\left(\begin{array}{ccc}\frac{\partial M_{j}(\mathbf{s}_i)}{\partial s}  \mathrm{v}_{1x}(\mathbf{s}_j) \frac{1}{J_s} (\pt{s}_i) -   M_{j} \mathrm{v}_{3x}\kappa_{x^{\prime}}& \frac{\partial M_{j}}{\partial s}  \mathrm{v}_{1y} \frac{1}{J_s} -   M_{j} \mathrm{v}_{3x}\kappa_{x^{\prime}}& \frac{\partial M_{j}}{\partial s}  \mathrm{v}_{1z}\frac{1}{J_s} -   M_{j} \mathrm{v}_{3z}\kappa_{x^{\prime}}\\ 
\\ \frac{\partial M_{j}}{\partial t}  \mathrm{v}_{2x} \frac{1}{J_t} -  M_{j} \mathrm{v}_{3x}\kappa_{y^{\prime}}& \frac{\partial M_{j}}{\partial t} \mathrm{v}_{2y} \frac{1}{J_t}-  M_{j} \mathrm{v}_{3y}\kappa_{y^{\prime}}& \frac{\partial M_{j}}{\partial t} \mathrm{v}_{tz} \frac{1}{J_t} -  M_{j} \mathrm{v}_{3z}\kappa_{y^{\prime}}\\
\\  \frac{\partial M_{j}}{\partial r}  \mathrm{v}_{3x} \frac{1}{J_r} & \frac{\partial M_{j}}{\partial r} \mathrm{v}_{3y} \frac{1}{J_r}& \frac{\partial M_{j}}{\partial r} \mathrm{v}_{3z} \frac{1}{J_r} \\  
\\ \frac{\partial M_{j}}{\partial t}  \mathrm{v}_{1x} \frac{1}{J_t}+ \frac{\partial M_{j}}{\partial s}  \mathrm{v}_{2x} \frac{1}{J_s}& \frac{\partial M_{j}}{\partial t}  \mathrm{v}_{1y} \frac{1}{J_t}+ \frac{\partial M_{j}}{\partial s}  \mathrm{v}_{2y}  \frac{1}{J_s}  &\frac{\partial M_{j}}{\partial t}  \mathrm{v}_{1z}  \frac{1}{J_{t}} + \frac{\partial M_{i}}{\partial s}  \mathrm{v}_{2z}  \frac{1}{J_s}  \\ \\ 0 & 0 & 0 \ \\ 0 & 0 & 0 \end{array}\right)
\nonumber
\end{equation}
In matrix notation we have:
\begin{equation}
\label{Straineq1}
 \myVecGreek{\epsilon}_{i}^{\prime}= \left[\begin{array}{ccc}\mathbf{ \hat{B}}^{\prime}_{i1} & \mathbf{ \hat{B}}^{\prime}_{i2}  & \cdots \end{array}\right] \left\{\begin{array}{c}\mathbf{ u}_{1} \\\mathbf{ u}_{2}  \\\vdots\end{array}\right\}
\end{equation}

\subsection{Use of NURBS for interpolation.}
If we use NURBS for interpolation functions $M_{j}$ we have to consider that $\mathbf{ u}_{j}$ are real values, whereas NURBS work with parameters.
This means that real displacements have to be converted to parameter values.
For the conversion we determine the NURBS parameters $\mathbf{ c}_{j}$ in such a way that the the real displacement values $\mathbf{u}(\pt{s}_{i})$ at all internal points locations $\pt{s}_{i}$ are exactly replicated:
\begin{equation}
\label{ }
\mathbf{ u}(\pt{s}_{i})= \sum_{j=1}^{J} M_{j}(\pt{s}_i) \mathbf{ c}_{j}
\end{equation}

If we gather all real values of displacements in vector $\{\myVec{u}\}$ the relationship between parameter values $\{\myVec{c}\}$ at these points and real values  is given by:
\begin{equation}
\label{ }
 \{\myVec{u}\} = [\mathbf{ A}]\{\myVec{c}\}
\end{equation}
The square conversion matrix $ [\mathbf{ A}]$ is given by:
\begin{equation}
\label{ }
 [\mathbf{ A}]= \left(\begin{array}{ccccccc}M_{1}(\pt{s}_{1})& 0 & 0 & M_{2}(\pt{s}_{1})& 0 & 0 & \cdots\\0 & M_{1}(\pt{s}_{1}) & 0 & 0 & M_{2}(\pt{s}_{1}) & 0 & \cdots\\0 & 0 & M_{1}(\pt{s}_{1}) & 0 & 0 & M_{2}(\pt{s}_{1}) & \cdots \\ M_{1}(\pt{s}_{2})& 0 & 0 & M_{2}(\pt{s}_{2})& 0 & 0 & \cdots\\0 & M_{1}(\pt{s}_{2}) & 0 & 0 & M_{2}(\pt{s}_{2}) & 0 & \cdots\\0 & 0 & M_{1}(\pt{s}_{2}) & 0 & 0 & M_{2}(\pt{s}_{2}) & \cdots \\ \vdots& \vdots & \vdots & \vdots& \vdots & \vdots & \ddots\end{array}\right)
\end{equation}
where $M_{i}(\pt{s}_{j})$ are NURBS basis function values.

The inverse relationship is given by:
\begin{equation}
\label{ }
\{\myVec{c}\}=  [\mathbf{ A}]^{-1}\{\myVec{u}\}
\end{equation}

The strains are now computed by:
\begin{equation}
\label{Straineq2}
 \myVecGreek{\epsilon}_{i}^{\prime}= \left[\begin{array}{ccc}\mathbf{ \hat{B}}^{\prime}_{i1} & \mathbf{ \hat{B}}^{\prime}_{i2}  & \cdots \end{array}\right] [\mathbf{ A}]^{-1}\left\{\begin{array}{c}\mathbf{ u}_{1} \\\mathbf{ u}_{2}  \\\vdots\end{array}\right\}
\end{equation}

\section{Test examples}
 The first two examples test the theory on plane inclusions which are loaded in perpendicular and tangential directions.  For one example we are able to compare with a closed form solution, for the other one we compare with a FEM result.
The examples relate to a bi-material cube of dimensions $1 \times 1 \times 1$ fixed at the bottom and subjected to a vertical tensile traction of 1 (all properties are non-dimensional).
The cube consists of 2 materials with Youngs moduli $E$ and $E_1$. 
For the first example we vary the thickness $d$ of the inclusion material as well as the ratio $E_1/E$. 
We compare the IGABEM results using options 1 (purely numerical integration) and 2 (combined numerical/analytical integration) with the analytical solution.
The third test example tests the ability to model curved surfaces. Here we can also compare with an analytical solution.
 For all runs $E=10$ was assumed. Poisson's ratio is assumed zero.

\subsection{Test example 1}
\begin{figure}
\begin{center}
\begin{overpic}[scale=0.4]{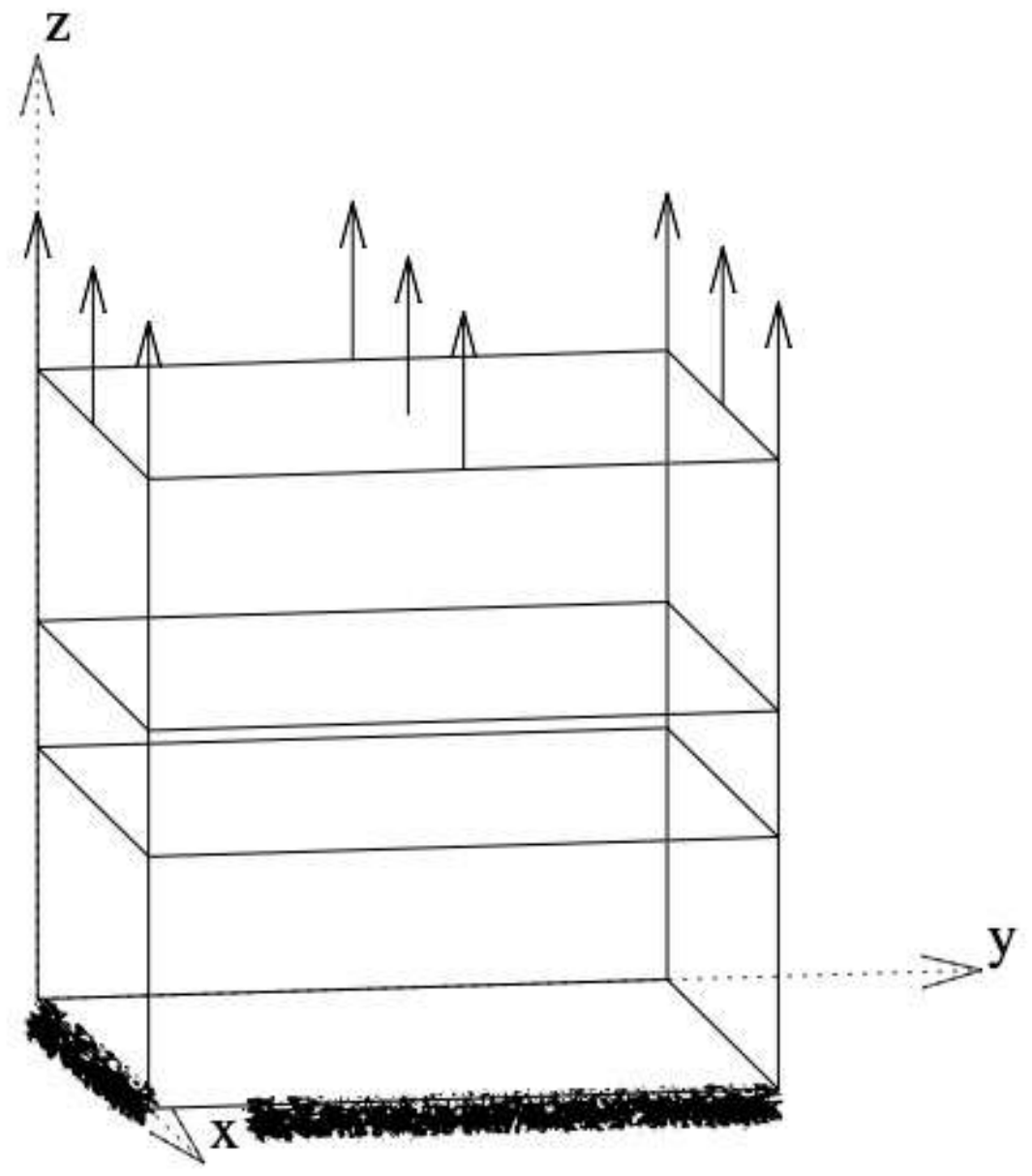} 
  \put(5,40){1}    
  \put(40,20){1}
  \put(65,13){1}
 \put(75,35){$d$}
\put(40,35){$E_1$}
\put(40,50){$E$}
\put(40,25){$E$}
\put(40,80){$t=1$}
\end{overpic}
\caption{Test example 1: Bi-material cube subjected to tensile load. Loading perpendicular to inclusion.}
\label{Test1}
\end{center}
\end{figure}

Here the inclusion is located perpendicular to the loading (see \myfigref{Test1}).
The theoretical value for the vertical displacement at the top, $u_z$,  is given by:
\begin{equation}
\label{ }
Eu_{z}=(E/E1 - 1)d +1
\end{equation}

\subsubsection{IGABEM discretisation}

\begin{figure}
\begin{center}
\begin{overpic}[scale=0.8]{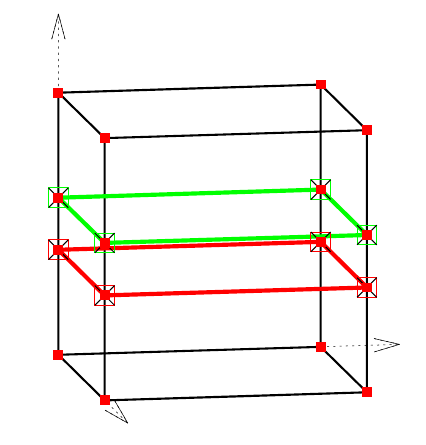} 
 \end{overpic}
\caption{Test example 1: IGABEM discretisation: Collocation points are shown as red squares, control points for the inclusion as colour coded hollow squares, inclusion surfaces depicted as coloured lines and inclusion points as crosses. }
\label{Test1IGA}
\end{center}
\end{figure}

The IGABEM discretisation is shown in \myfigref{Test1IGA}.
It consists of 6 linear surface patches. The inclusion is specified by 2 linear surfaces. To change the continuity of displacement at the points where the inclusion meets the boundary to $C^0$, knots were inserted into the basis functions describing the geometry to arrive at the basis functions for approximating the unknown, resulting in the collocation points shown.

\subsubsection{Results}
We compare in table \ref{tab:Test1} the IGABEM results using options 1 and 2 for different inclusion thickness and ratio of moduli.
\begin{mytable}
  {H}               
  {Results of test 1:vertical displacement at top}  
  {tab:Test1}  
  {ccccc}         
  \mytableheader{ $E_1/E$ &  d & IGABEM option 1& IGABEM option 2 & theory }  
2 & 0.2 & 0.09 & 0.09 & 0.09 \\
4 & 0.1 & 0.0925 & 0.0925 & 0.0925 \\
8 & 0.05 & 0.09559 & 0.0956  & 0.09562 \\
8 & 0.025 & 0.1 & 0.0978  & 0.0978 \\
\end{mytable}%
Since for this example the initial stress is constant inside the inclusion, the results for option 1 and 2 are identical for larger thicknesses. However, if the thickness is very small the results of option 1 are in error because the bad aspect ratio of the tetrahedral subregions.
\begin{figure}[H]
\begin{center}
\begin{overpic}[scale=0.7]{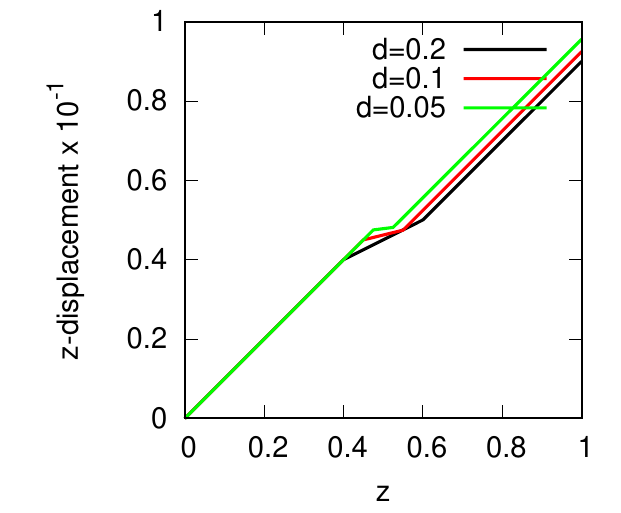} 
 \end{overpic}
\caption{Test example 1: Vertical displacement along z-axis }
\label{Displ1}
\end{center}
\end{figure}
A plot of the vertical displacements in a vertical direction is shown in \myfigref{Displ1}.

\newpage

\subsubsection{Test example 2}
\begin{figure}
\begin{center}
\begin{overpic}[scale=0.4]{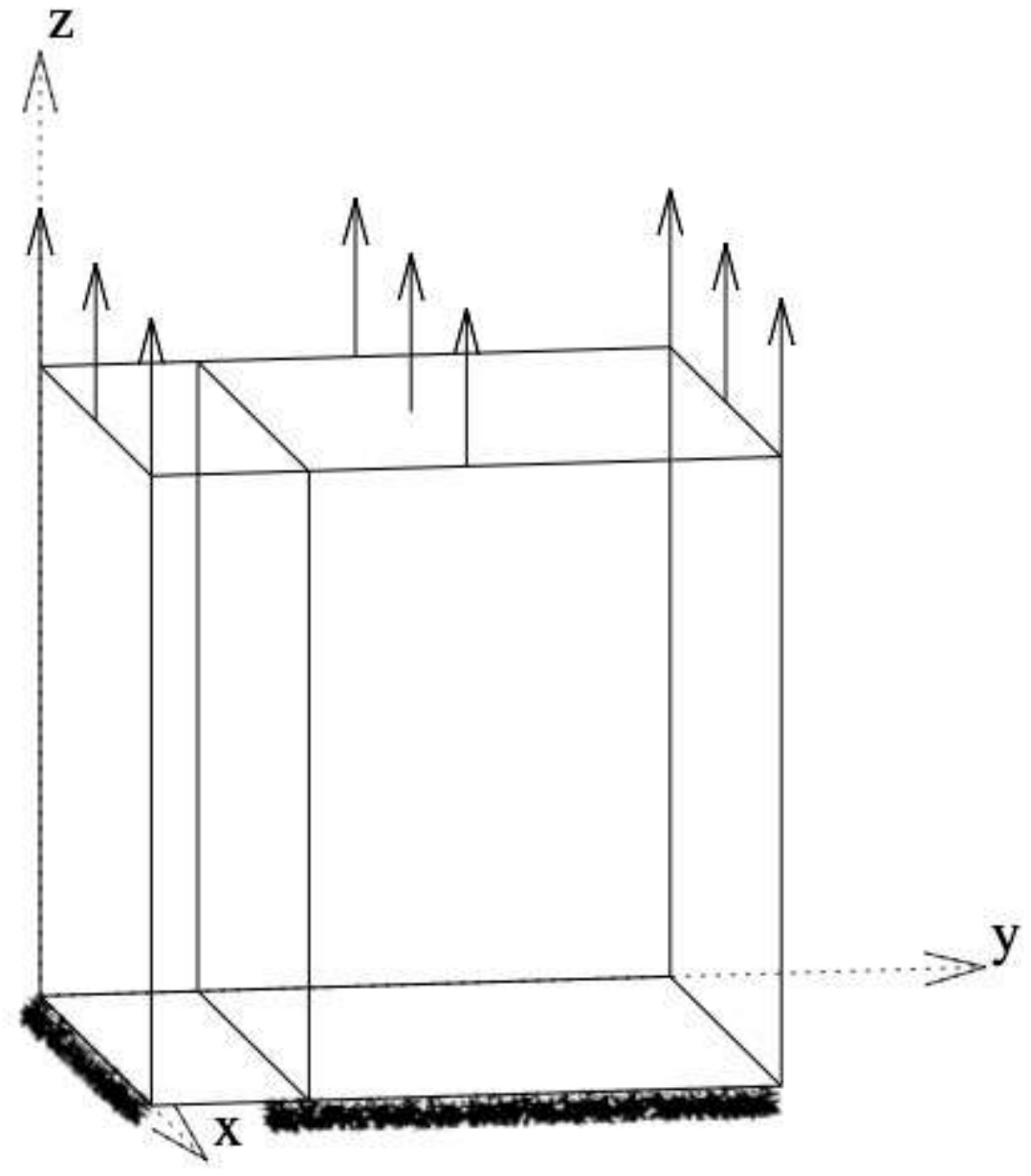} 
  \put(5,40){1}    
  \put(40,20){1}
  \put(65,13){1}
 \put(25,52){$d$}
\put(15,40){$E_1$}
\put(40,40){$E$}
\put(40,80){$t=1$}
\end{overpic}
\caption{Test example 2: Bi-material cube subjected to tensile load. Loading parallel to inclusion.}
\label{Test2}
\end{center}
\end{figure}

To test that the method works also for less simple examples now the inclusion is located tangential to the loading (see \myfigref{Test2}).
There is no theoretical solution so comparison is done with results of Finite Element analyses.
\subsubsection{IGABEM discretisation}

\begin{figure}
\begin{center}
\begin{overpic}[scale=0.8]{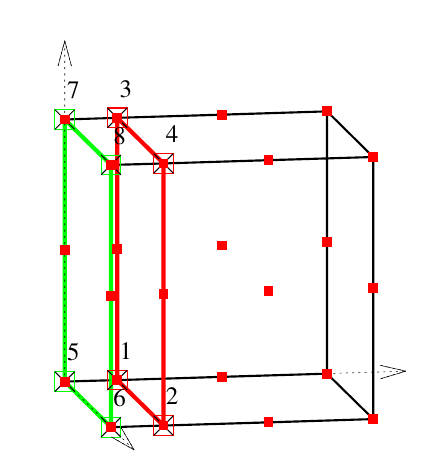} 
 \end{overpic}
\caption{Test example 2: IGABEM discretisation: Collocation points are shown as red squares, control points for the inclusion as colour coded hollow squares, inclusion surfaces with coloured lines and inclusion points as crosses. }
\label{Test2IGA}
\end{center}
\end{figure}

The IGABEM discretisation is shown in \myfigref{Test2IGA}.
It consists of 10 linear surface patches. The inclusion is specified by 2 linear surfaces. To determine the basis functions for the approximation of the unknown displacements we order elevate from linear to quadratic in 2 directions, resulting in the collocation points shown. The model has 72 degrees of freedom.

\subsubsection{FEM discretisation}
\begin{figure}
	\begin{center}
		\begin{overpic}[scale=0.7]{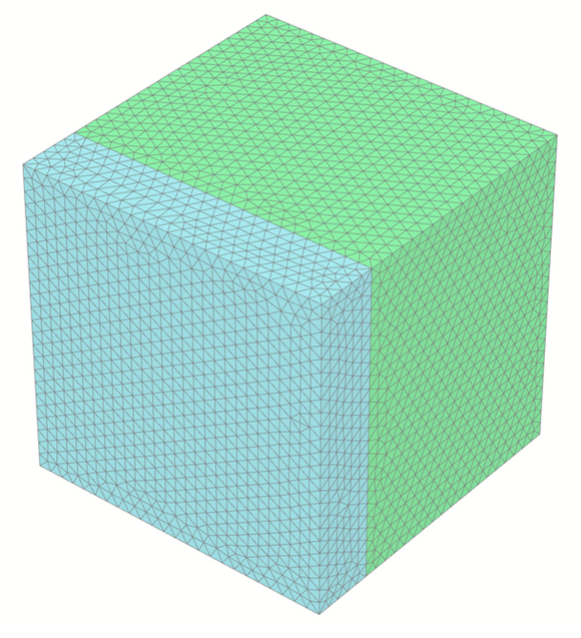} 
		\end{overpic}
		\caption{Test example 2: FEM discretisation }
		\label{FEM_Mesh}
	\end{center}
\end{figure}

The FEM mesh is shown in \myfigref{FEM_Mesh}. It consits of 10-noded quadratic tetrahedral elements. The model has approximately 270.000 degrees of freedom.

\newpage

\subsubsection{Results}
\begin{figure}
\begin{center}
\begin{overpic}[scale=0.7]{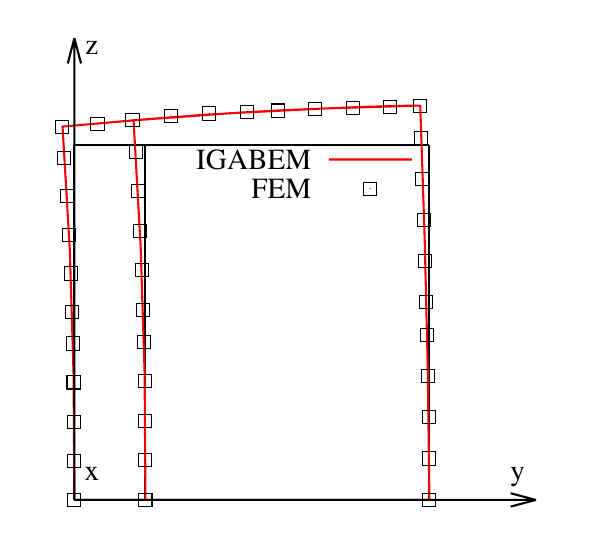} 
 \end{overpic}
\caption{Test example 2: Displaced shape for $E_1/E=2$ and $d=0.2$. Comparison between IGABEM and FEM}
\label{Test2Displ}
\end{center}
\end{figure}

The displaced shape for $E_1/E=2$ and $d=0.2$ is shown in \myfigref{Test2Displ}. Close agreement between the IGABEM and FEM results can be observed.

\subsection{Test Example 3}
The example is an infinite cylindrical opening in an infinite domain subjected to an hydrostatic pressure of $p_i=1$.
The radius $R$ of the opening is 5. 
The domain is reinforced by a cylindrical lining with a thickness $d=0.5$ and $E_1/E=2$.
This example is designed to show the advantage of using local strains and interpolation with NURBS. 
\begin{figure}[h]
\begin{center}
\begin{overpic}[scale=0.11]{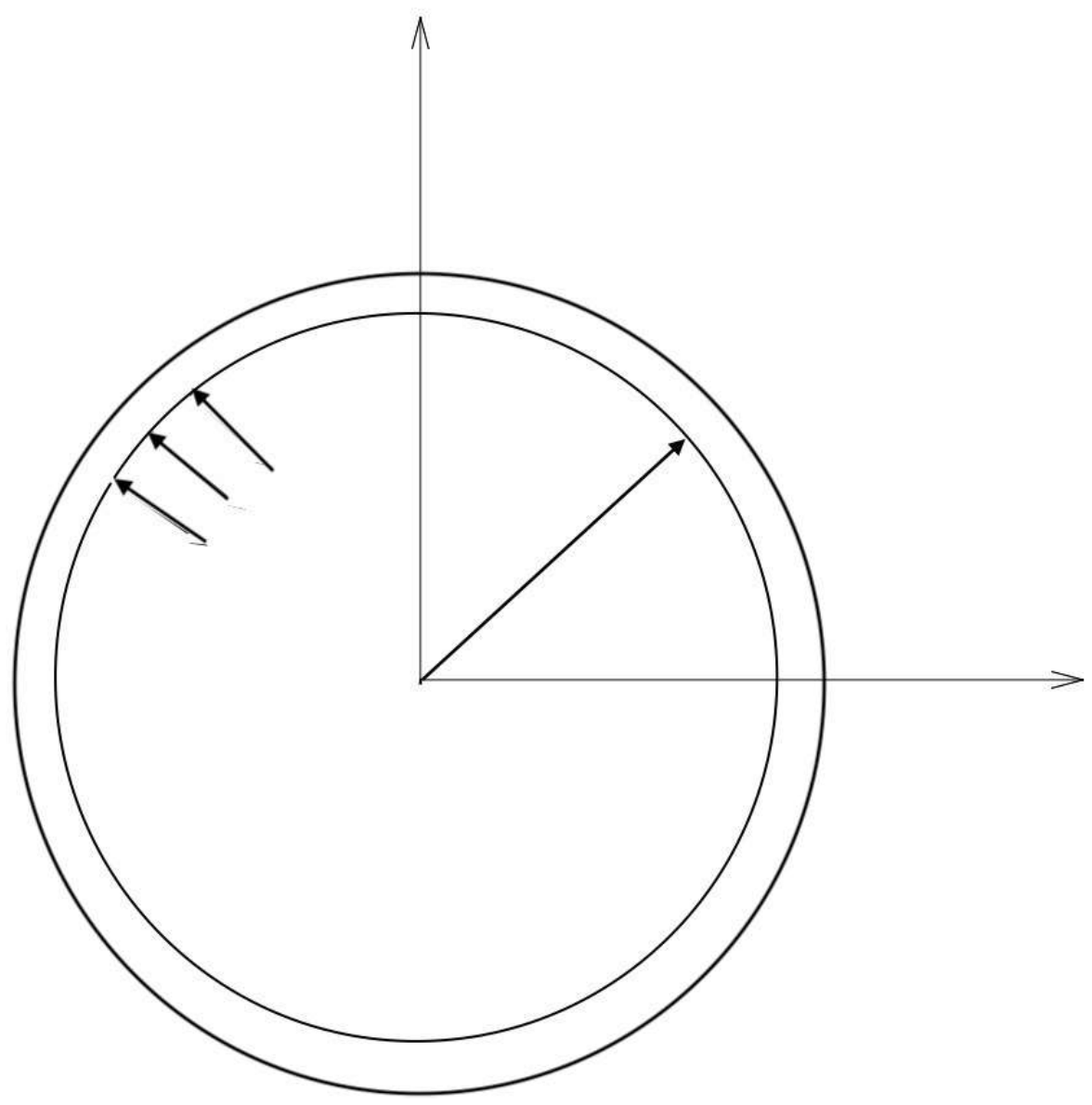}
 \put(60,45){$R$}
\end{overpic}
\begin{overpic}[scale=0.25]{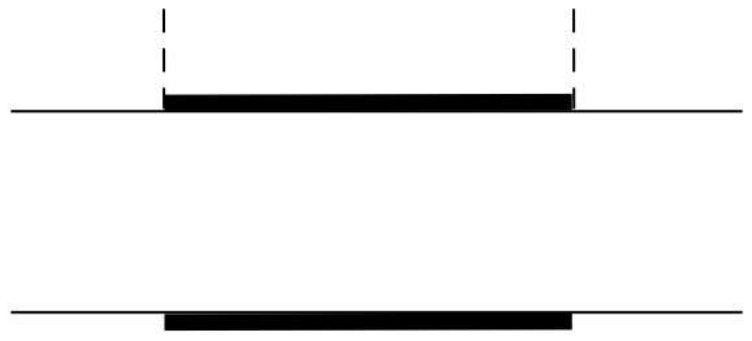}
    \put(10,20){plane strain}
    \put(65,20){plane strain}
\end{overpic}
\caption{Cross-section and side view of test example.}
\label{TestShot}
\end{center}
\end{figure}

For the 3-D analysis we assume that the extent of the reinforcement ring along axis of the cylindrical opening is $2R$. Plane strain conditions (i.e. the assumption is that displacements remain constant along the axis of the opening) are assumed where the reinforcement finishes.

\subsubsection{Analytical solution}
The numerical results can be compared with a plane strain analytical solution. To obtain the analytical solution of the bi-material problem, we use two solutions: One for the circular hole in an infinite domain with the modulus of $E$ and one for the thick cylinder with the modulus $E_1$. $R$ is the radius of the circular hole and $R_0$ the radius to the interface. Therefore $d=R_o - R$.
Plane strain conditions and $\nu=0$ are assumed.
\begin{figure}[h]
	\begin{center}
		\begin{overpic}[scale=0.15]{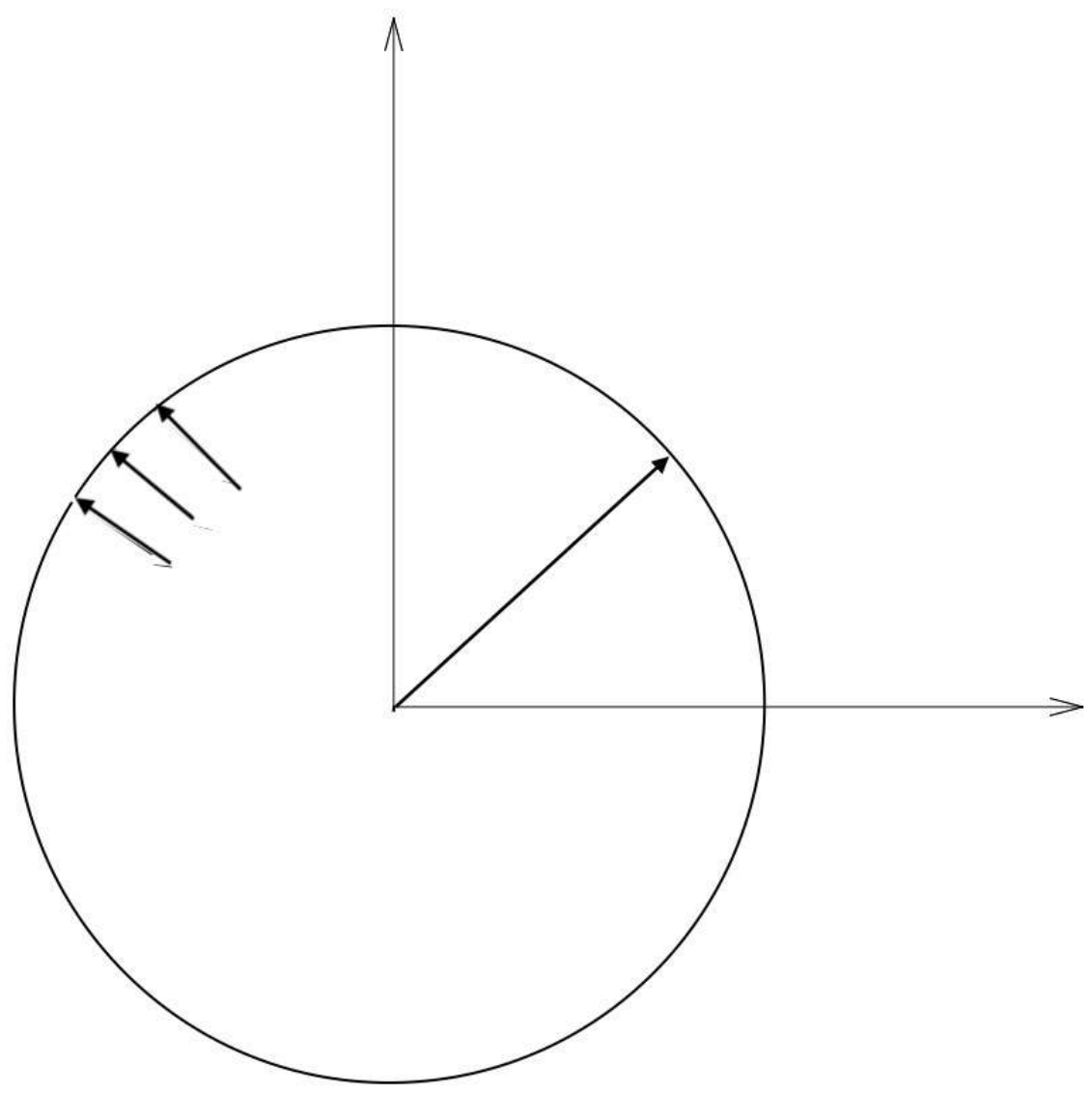}		
			\put(25,60){$p_{oK}$}
			\put(60,50){$R_o$}
			\put(50,70){$E$}
		\end{overpic}
		\begin{overpic}[scale=0.14]{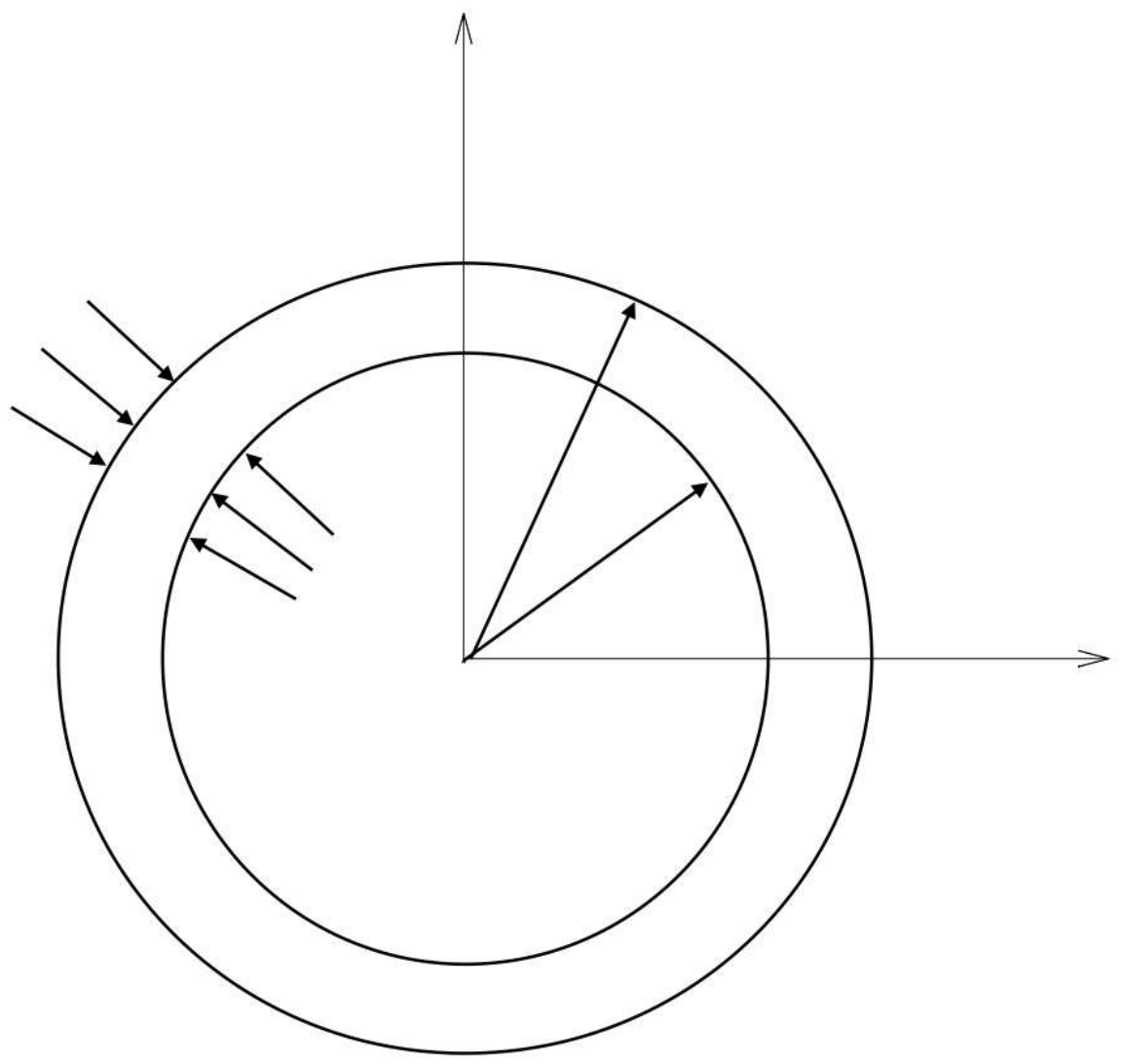}		
			\put(30,50){$p_i$}
			\put(18,65){$p_{oT}$}
			\put(50,45){$R$}
			\put(45,55){$R_o$}
			\put(10,30){$E_1$}
		\end{overpic}
		\caption{For getting the exact solution of the problem we combine the circular excavation (Kirsch) problem (left) with the thick cylinder problem (right).}
		\label{TestShot}
	\end{center}
\end{figure}

The radial displacement of the circular boundary at $R_o$, due to a hydrostatic pressure of $p_{oK}$, is according to Kirsch \cite{Kirsch}:
\begin{equation}
	\label{ }
	u_{rK}= \frac{R_o~p_{oK}}{E}
\end{equation}
The radial displacement at $R_o$ of the thick cylinder due to $p_{oT}$  and $p_i$ is:
\begin{equation}
	\label{ }
	u_{rT}= \frac{1}{E_1}~\frac{R_o}{R_o^2~-~R^2}\left(2~R^2~p_i~-~\left(R^2~+~R_o^2\right)~p_{oT}\right) 
\end{equation}

Equilibrium at the interface 
\begin{equation}
	\label{ }
	p_{oK}=p_{oT}=p_{o}
\end{equation}
and compatibility
\begin{equation}
	\label{ }
	u_{rT}=u_{rK}
\end{equation}
gives
\begin{equation}
	\label{ }
	\frac{R_o~p_{o}}{E}~=~\frac{1}{E_1}~\frac{R_o}{R_o^2~-~R^2}\left(2~R^2~p_i~-~\left(R^2~+~R_o^2\right)~p_{o}\right) 
\end{equation}
which can be solved for $p_o$:
\begin{equation}
	\label{ }
	p_o~=~\frac{2~E~R^2~p_i}{E~\left(R^2~+~R_o^2~\right)~+~E_1~\left(R_o^2-R^2\right)    }
\end{equation}

The radial displacement inside the thick cylinder at  a distance $r$ is obtained by:
\begin{equation}
	\label{ }
	u_{rT}= \frac{R^2~p_i~\left[E~\left(R_o^2~-~r^2\right)~+~E_1~\left(R_o^2~+~r^2\right)\right]}{E_1~r~\left[E~\left(R_o^2~+~R^2\right)~+~E_1~\left(R_o^2~-~R^2\right)\right]}
\end{equation}
and inside the infinite domain by:
\begin{equation}
	\label{ }
	u_{rK}= \frac{p_o~R_o^2}{E~r}
\end{equation}

\subsubsection{IGABEM model}

\begin{figure}
\begin{center}
\includegraphics[scale=0.7]{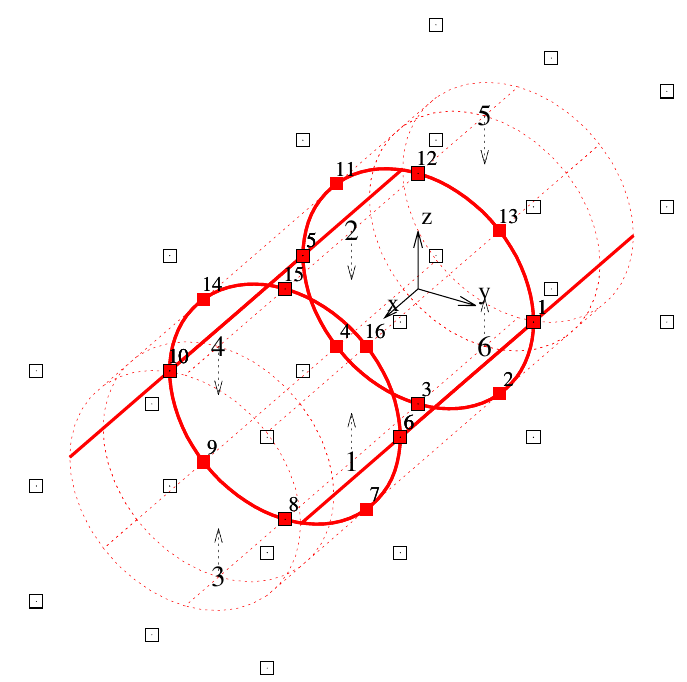}
\caption{Test example 3: Surface discretisation with finite and infinite patches. Control points are shown as hollow squares. Collocation points are shown as red filled squares. Also shown are the outward normal vectors and the numbering of the collocation points and patches. }
\label{Smesh}
\end{center}
\end{figure}

The surface discretisation is shown in \myfigref{Smesh}. It consists of 2 finite and 4 infinite (plane strain) patches. The NURBS surfaces are of order 2 in the circumferential direction and of order 1 in the direction along the opening. The model has 48 degrees of freedom. The geometry of the cylindrical opening is exactly replicated.

\begin{figure}[H]
\begin{center}
\includegraphics[scale=0.6]{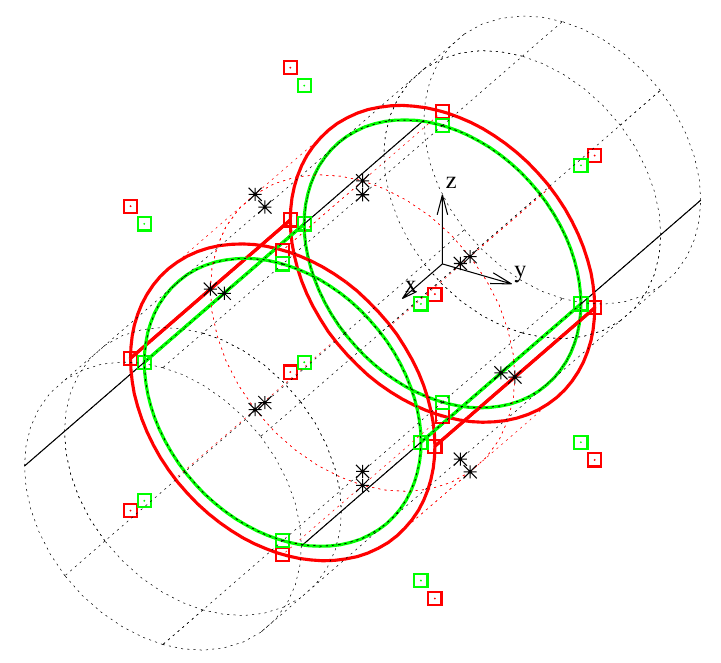}
\caption{Test example 3: Definition of reinforcement ring. The control points defining the inclusions are depicted by coloured red and green squares and the bounding surfaces by colored lines. Internal points are depicted by stars. }
\label{Incl}
\end{center}
\end{figure}

The definition of the reinforcement ring is depicted in \myfigref{Incl}. 
It consists on 2 inclusions defined by NURBS surfaces of order 2 in the circumferential direction and of order 1 in the direction along the opening. The geometry of the ring including its curvature is exactly replicated. Only 32 internal points are defined allowing linear variation of the displacements in radial direction and constant variation in the direction along the axis of the opening.

\subsubsection{Results}
\begin{figure}[H]
\begin{center}
\includegraphics[scale=0.6]{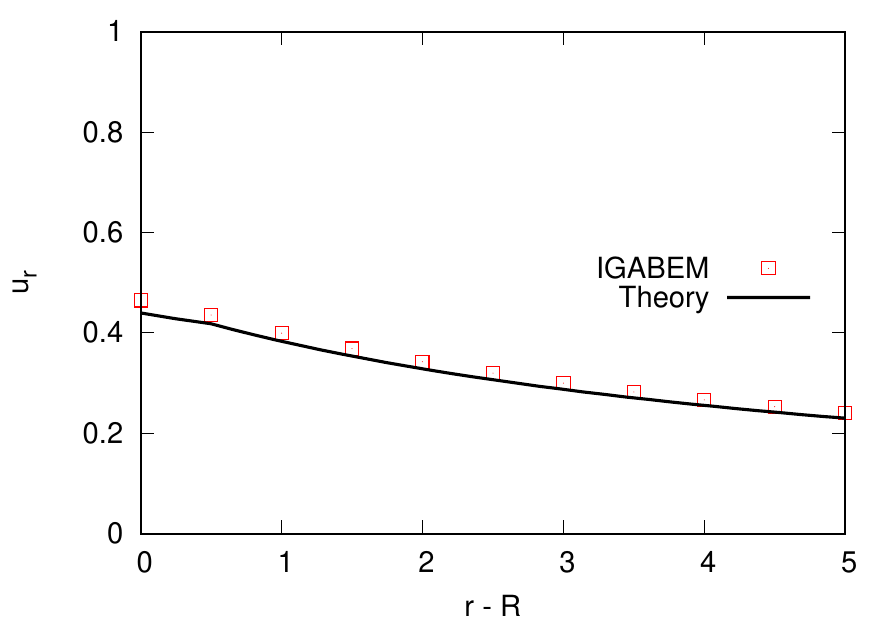}
\caption{Test example 3: Comparison of the variation of radial displacement in a radial direction from the surface of the opening. }
\label{Res}
\end{center}
\end{figure}

For this example the circumferential strain inside the inclusion does not vary along the circumference. Using NURBS of order 2 for the interpolation of the displacements inside the inclusion and considering  that the curvature is exactly computed, the variation of the strains in circumferential direction is therefore also exactly computed. Errors in the results are only due to the fact that we assumed a linear variation of the displacements across the thickness which results in a constant stress (the theoretical solution gives a linearly varying stress) and due to the numerical integration.
In \myfigref{Res} we compare the radial displacement obtained by the simulation with the theory. Good agreement can be found.

\section{Summary and Conclusions}
The simulation with the isogeometric BEM of problems that contain inclusions has been discussed. In order to simulate such problems with the BEM the original boundary integral equations have to be amended. An efficient and accurate method is to treat the inclusions as subdomains, where initial stresses are generated due to the difference in material properties. The authors have published on this topic previously and have demonstrated that the method works well as long as the inclusions have a moderate aspect ratio. If one of the dimensions of the inclusion is much smaller than the other dimension (we referred to these as thin inclusions) the methods do not give good results and in extreme cases may give wrong answers.
This is because the published integration methods break down for extreme aspect ratios. 

In this paper we have presented efficient and accurate methods for thin inclusions.
This includes a robust and efficient integration scheme for the singular integrals that occur. It is also demonstrated that using a local strain definition substantial improvements in efficiency and accuracy can be achieved for inclusions which are curved. The formula for computing the strain on a surface from displacement values, which has appeared in a few publications, is amended to account for the influence of curvature.

Three test examples are presented that attest the robustness, accuracy and efficiency of the implementation.
Applications are in the area of reinforced structures or concrete linings in underground construction. Given our expertise in this area the latter will be the main area of application.

Finally it should be noted that the method of initial stresses used is ideally suited for implementing non-linear (elasto-plastic) material behaviour. All that is is required is a different formula for computing the initial stress. This has already been implemented for normal inclusions in \cite{BEER2021} and can be applied to thin inclusions without any changes.

\appendix
\section{Analytical singular solutions}

 \begin{eqnarray}
 \label{eq_DeltaEprime_11}
 \triangle \fund{E^{\prime}}(1,1) & = &
\frac{1}{48} C \left\{ R \cos{3 \tilde{\theta}} \left( \sin{3 \phi_2} -
    9 \sin{\phi_1} - \sin{3 \phi_1} \right) + \sin{\phi_2} \left[ 24 C_L \, H \left( 3 +  \right. \right. \right. \\
    & & \left. \left. \left. + 4 C_3 + \cos{2 \phi_2} \right) -
 6 R (7 + 16 C_3 + 3 \cos{2 \phi_2}) \cos{\tilde{\theta}} + 9 R \cos{3 \tilde{\theta}} + \right. \right. \nonumber \\
 & & \left. \left. + 6 H (17 + 16 C_3 + 5 \cos{2 \phi_2}) \sin{\tilde{\theta}} - 2 H (5 + \cos{2\phi_2}) 
    \sin{3 \tilde{\theta}} \right] + \right. \nonumber \\
    & & \left. 2 \sin{\phi_1} \left[ -12 C_L H (3 + 4 C_3 + 
    \cos{2 \phi_1}) + 3 R (7 + 16 C_3 + 3 \cos{2 \phi_1})
     \cos{\tilde{\theta}} - \right. \right. \nonumber \\ 
     & & \left. \left. 3 H (17 + 16 C_3 + 5 \cos{2 \phi_1})
  \sin{\tilde{\theta}} + H (5 + \cos{2 \phi_1}) \sin{3 \tilde{\theta}} \right] \right\}
\nonumber \\
 \label{eq_DeltaEprime_12}
 \triangle \fund{E^{\prime}}(1,2) & = & 
\frac{1}{48} C \left\{ 4 R \left[ \cos{3 \tilde{\theta}} \left( \sin^3{\phi_2} - \sin^3{\phi_1} \right) + 3 \cos{\tilde{\theta}} \left( 4 \sin{\phi_2} - 3 \sin^3{\phi_2} - 4 \sin{\phi_1} +   \right. \right. \right. \\
    & & \left. \left. \left. + 3 \sin^3{\phi_1} \right) \right] + H \left[ - 12 C_L \left( \sin{\phi_2} + \sin{3 \phi_2} - 4 \sin{\phi_1} \cos^2{\phi_1} \right) + 3 \left( \sin{\phi_1} + \right. \right. \right. \nonumber \\
& & \left. \left. \left. + 5 \sin{3 \phi_1} - \sin{\phi_2} - 5 \sin{3 \phi_2} \right) \sin{\tilde{\theta}} + 4 \left(  \sin^3{\phi_1} - \sin^3{\phi_2} \right) \sin{3 \tilde{\theta}}  \right] \right\} \nonumber \\
 \label{eq_DeltaEprime_13}
 \triangle \fund{E^{\prime}}(1,3) & = &
- C \left( \sin{\phi_2} - \sin{\phi_1} \right) \left[ C_L H + \right.
\\
& &  \left. \sin{\tilde{\theta}} \left( H - R \cos{\tilde{\theta}} \sin{\tilde{\theta}} + H \sin^2{\tilde{\theta}} \right) \right]
\nonumber \\
 \label{eq_DeltaEprime_14}
 \triangle \fund{E^{\prime}}(1,4) & = &
\frac{1}{24} C \left\{ 4 R \left[ 3 \left( 4 C_3 \cos{\phi_2} + 3 \cos^3{\phi_2} - 4 C_3 \cos{\phi_1} - 3 \cos^3{\phi_1} \right) \cos{\tilde{\theta}} \right. \right. \\
& & \left. \left. + \left( \cos^3{\phi_1} - \cos^3{\phi_2} \right) \cos{3 \tilde{\theta}} \right] + H \left[ 6 \cos{\phi_2}  \left( -4 C_L \left( 1 +
    2 C_3 + \cos{2 \phi_2} \right) - \right. \right. \right. \nonumber \\ 
    & & \left. \left. \left. - \left( 5 + 8 C_3 + 5 \cos{2 \phi_2} \right)  \sin{\tilde{\theta}} \right] + 6  \cos{\phi_1} \left[ 8 C_L \left( C3+cos^2{\phi_1} \right) + \right. \right. \right. \nonumber \\
& & \left. \left. \left.    \left( 5 + 8 C_3 + 5 \cos{2 \phi_1} \right) \sin{\tilde{\theta}} \right] + 4 \left( \cos^3{\phi_2} - \cos^3{\phi_1} \right) \sin{3 \tilde{\theta}} \right] \right\};
\nonumber \\
 \label{eq_DeltaEprime_15}
 \triangle \fund{E^{\prime}}(1,5) & = &
\frac{1}{8} C \left( \cos{2 \phi_2} - \cos{2 \phi_1} \right) \left( -8 H - 4 R + 9 H \cos{\tilde{\theta}} - \right. \\
& &  \left. - H \cos{3 \tilde{\theta}} + 4 R \sin^3{\tilde{\theta}} \right)
\nonumber \\
 \label{eq_DeltaEprime_16}
 \triangle \fund{E^{\prime}}(1,6) & = &
\frac{1}{8} C \left\{ H \left[ 4 \left( \phi_2 - \phi_1 \right) \cos{\tilde{\theta}} \left( \cos{2 \tilde{\theta}} - 5 - 4 C_3 \right) + \left( \cos{3 \tilde{\theta}} - 9 \cos{\tilde{\theta}} \right) \sin{2 \phi_2} + \right. \right. \\
& & \left. \left. + 8 \left( 2 (1+C_3) \left( \phi_2 - \phi_1 \right) + \sin{2 \phi_2} - \sin{2 \phi_1} \right) - 2 \cos{\tilde{\theta}} \left( \cos{2 \tilde{\theta}} - 5 \right) \sin{2 \phi_1} \right] - \right. \nonumber \\
& & \left. -4 R 
\left[ \left( \phi_2 - \phi_1 \right) \left( \sin{\tilde{\theta}} - 1 \right) \left( 3 + 4 C_3 - \cos{2 \tilde{\theta}} + 2 \sin{\tilde{\theta}} \right) + \sin{2 \phi_2} \left(  \sin^3{\tilde{\theta}} - 1 \right) - \right. \right. \nonumber \\
& & \left. \left. - \sin{2 \phi_1} \left( \sin^3{\tilde{\theta}} - 1 \right) \right] \right\}
\nonumber \\
 \label{eq_DeltaEprime_21}
 \triangle \fund{E^{\prime}}(2,1) & = &
\frac{1}{48} C \left\{  3 R \left( 3 \cos{3\phi_2} - 7 \cos{\phi_2} + 7 \cos{\phi_1} - 3 \cos{\phi_1} \right) \cos{\tilde{\theta}} + 4 R \left( \cos^3{\phi_1} - \right. \right. \\
& & \left. \left. - \cos^3{\phi_2} \right) \cos{3 \tilde{\theta}} + 12 C_L H \left( \cos{3 \phi_1} - \cos{\phi_1} + 4 \cos{\phi_2} \sin^2{\phi_2}  \right) + 3 H \left( \cos{\phi_2}  \right. \right. \nonumber \\
& & \left. \left. - 5 \cos{3 \phi_2} - \cos{\phi_1} + 5 \cos{3 \phi_1} \right) \sin{\tilde{\theta}} + 4 H \left( \cos^3{\phi_2} - \cos^3{\phi_1}  \right) \sin{3 \tilde{\theta}} \right\}
\nonumber \\
 \label{eq_DeltaEprime_22}
 \triangle \fund{E^{\prime}}(2,2) & = &
\frac{1}{48} C \left\{ - 12 C_L H \left[  \cos{\phi_2} \left( 6 + 8 C_3 - 2 \cos{2 \phi_2} \right) + 2 \cos{\phi_1} \left( \cos{2 \phi_1} - 3 - 4 C_3   \right) \right] - \right. \\
& &  \left. - 6 R \left( \cos{\phi_2} - \cos{\phi_1}  \right) \left(  3 \cos{2 \phi_2} - 4 - 16 C_3 + 6 \cos{\phi_2} \cos{\phi_1} + 3 \cos{2 \phi_1} \right) \cos{\tilde{\theta}} + \right. \nonumber \\
& & \left. + R \left( \cos{3 \phi_2} - 9 \cos{\phi_2} + 9 \cos{\phi_1}  - \cos{3 \phi_1} \right) \cos{3 \tilde{\theta}} + 6 H \left( \cos{\phi_2} - \cos{\phi_1}  \right) \left( 5 \cos{2 \phi_2} - 12 - 16 C_3 + \right. \right. \nonumber \\
& & \left. \left.  + 10 \cos{\phi_2} \cos{\phi_1} + 5 \cos{2 \phi_1}   \right) \sin{\tilde{\theta}} + H \left(  9 \cos{\phi_2} - \cos{3 \phi_2} - 9 \cos{\phi_1} + \cos{3 \phi_1}  \right) \sin{3 \tilde{\theta}} \right\}
\nonumber \\
 \label{eq_DeltaEprime_23}
 \triangle \fund{E^{\prime}}(2,3) & = & 
 C \left( \cos{\phi_2} - \cos{\phi_1} \right)  \left[  C_L H + \sin{\tilde{\theta}} \left(  H - R \cos{\tilde{\theta}} \sin{\tilde{\theta}} + H \sin^2{\tilde{\theta}} \right) \right]
\nonumber \\
 \label{eq_DeltaEprime_24}
 \triangle \fund{E^{\prime}}(2,4) & = & 
\frac{1}{24} C \left\{  4 R \left[  \cos{3 \tilde{\theta}} \left( \sin^3{\phi_2} - 
\sin^3{\phi_1} \right) + 3 \cos{\tilde{\theta}} \left( - 4 C_3 \sin{\phi_2}
- 3 \sin^3{\phi_2} +  \right. \right. \right. \\
& & \left. \left. \left. + 4 C_3 \sin{\phi_1} + 3 \sin^3{\phi_1}   \right) \right] + H 
\left[  6 \sin{\phi_2} \left( - 4 C_L \left( - 1 - 2 C_3 + \cos{2 \phi_2}    
\right) + \left(  5 + 8 C_3 - \right. \right. \right. \right. \nonumber \\
&  & \left. \left. \left. \left. - 5 \cos{2 \phi_2}  \right)
 \sin{\tilde{\theta}}  \right) + 6 \sin{\phi_1} \left( - 8 C_L \left(  C_3 + \sin^2{\phi_1}  \right) + \sin{\tilde{\theta}} \left( 5 \cos{2 \phi_1} - 5 - 8 C_3   \right) \right) \right. \right. \nonumber \\
 & & \left. \left. + 4 \left( \sin^3{\phi_1} - \sin^3{\phi_2}   \right) \sin{3 \tilde{\theta}} \right] \right\}
\nonumber \\
 \label{eq_DeltaEprime_25}
 \triangle \fund{E^{\prime}}(2,5) & = & 
\frac{1}{8} C \left\{ H \cos{\tilde{\theta}} \left[ 4 \left( 5 + 4 C_3 \right) \left(\phi_1 - \phi_2 \right) + 10 \sin{2 \phi_2} - 2 \cos{2 \tilde{\theta}} \left(- 2 \phi_2 + 2 \phi_1 + \right. \right. \right. \\
& & \left. \left. \left. \sin{2 \phi_2} \right) - 9 \sin{2 \phi_1} \right] + \sin{2 \phi_1} \left(8 H + 4 R + H \cos{3 \tilde{\theta}} - 4 R \sin^3{\tilde{\theta}} \right) - \right. \nonumber \\
& & \left. - 4 \sin{2 \phi_2} \left( 2 H + R - R \sin^3{\tilde{\theta}} \right) + 2 \left(\phi_2 - \phi_1 \right) \left[ 4 R + 8 \left(  H + C_3 H + C_3 R \right) - \right. \right. \nonumber \\ 
& & \left. \left. \left(  3 + 8 C_3 \right) R \sin{\tilde{\theta}} + R \sin{3 \tilde{\theta}}\right]
\right\}
\nonumber \\
 \label{eq_DeltaEprime_26}
 \triangle \fund{E^{\prime}}(2,6) & = & 
\frac{1}{8} C \left( \cos{2 \phi_2} - \cos{2 \phi_1}  \right) \left( 9 H \cos{\tilde{\theta}} - 8 H - 4 R  - H \cos{3 \tilde{\theta}} + 4 R \sin^3{\tilde{\theta}} \right)
\end{eqnarray}
\begin{eqnarray}
 \label{eq_DeltaEprime_31}
 \triangle \fund{E^{\prime}}(3,1) & = & 
\frac{1}{16} C \left[ 2 H \cos{\tilde{\theta}} \left( \phi_1 - \phi_2 - 9 \cos{\phi_2} \sin{\phi_2} + 9 \cos{\phi_1} \sin{\phi_1}   \right) + H \cos{3 \tilde{\theta}} \left( 2 \phi_2 - 2 \phi_1 + \right. \right. \\
& & \left. \left. + \sin{2 \phi_2} -  \sin{2 \phi_1} \right) + 4 \sin{2 \phi_2} \left( 2 H + R - R \sin^3{\tilde{\theta}}  \right) - 4 \sin{2 \phi_1} \left( 2 H + R - R \sin^3{\tilde{\theta}}  \right) + \right. \nonumber \\ 
& & \left. + 2 \left( \phi_2 - \phi_1   \right) R \left( 5 \sin{\tilde{\theta}} - 4 + \sin{3 \tilde{\theta}}   \right) \right]
\nonumber \\
 \label{eq_DeltaEprime_32}
 \triangle \fund{E^{\prime}}(3,2) & = & 
- \frac{1}{16} C \left[ 2 H \cos{\tilde{\theta}} \left( \phi_2 - \phi_1 - 9 \cos{\phi_2} \sin{\phi_2} + 9 \cos{\phi_1} \sin{\phi_1} \right) + \right. \\ 
& & \left. + H \cos{3 \tilde{\theta}} \left( 2 \phi_1 - 2 \phi_2 + \sin{2 \phi_2} - \sin{2 \phi_1}   \right) + 8 \left( \phi_2 - \phi_1  \right) R \left( 1 - 2 \sin{\tilde{\theta}} + \sin^3{\tilde{\theta}}   \right) + \right. \nonumber \\
& &  \left. + 4 \sin{2 \phi_2} \left( 2 H + R - R \sin^3{\tilde{\theta}}  \right) - 4 \sin{2 \phi_1} \left( 2 H + R - R \sin^3{\tilde{\theta}}   \right)  \right]
\nonumber \\
 \label{eq_DeltaEprime_33}
 \triangle \fund{E^{\prime}}(3,3) & = & 
C \left( \phi_1 - \phi_2   \right) H \left( \cos{\tilde{\theta}} - 1  \right) \left( 2 C_3 + \cos{\tilde{\theta}} + \cos^2{\tilde{\theta}}  \right) + \\
& &  + \frac{1}{2} C \left( \phi_2 - \phi_1   \right) R \left[ 2 + 4 C_3 - \left( 3 + 4 C_3 + \cos{2 \tilde{\theta}} \right) \sin{\tilde{\theta}}  \right]
\nonumber \\
\label{eq_DeltaEprime_34}
 \triangle \fund{E^{\prime}}(3,4) & = & 
\frac{1}{8} C \left( \cos{2 \phi_2} - \cos{2 \phi_1}  \right) \left( 9 H \cos{\tilde{\theta}} - 8 H - 4 R - H \cos{3 \tilde{\theta}} + 4 R \sin^3{\tilde{\theta}}  \right)\\
 \label{eq_DeltaEprime_35}
 \triangle \fund{E^{\prime}}(3,5) & = & 
2 C \left( \cos{\phi_2} - \cos{\phi_1}  \right) \left[ R \cos{\tilde{\theta}} \left( C_3 + \cos^2{\tilde{\theta}}  \right) + H \left( \sin^3{\tilde{\theta}} - C_3 C_L - C_3 \sin{\tilde{\theta}}  \right) \right]  \\
 \label{eq_DeltaEprime_36}
 \triangle \fund{E^{\prime}}(3,6) & = & 
- 2 C \left( \sin{\phi_2} - \sin{\phi_1}  \right) \left[ C_3 R \cos{\tilde{\theta}} + R \cos^3{\tilde{\theta}} + H \left( \sin^3{\tilde{\theta}} - \right. \right. \\
& & \left. \left. - C_3 C_L - C_3 \sin{\tilde{\theta}}  \right) \right]
\nonumber
\end{eqnarray}
where:
\begin{equation}
C_L = \ln{\left[\frac{\cos{\tilde{\theta}/2}-\sin{\tilde{\theta}/2}}{\cos{\tilde{\theta}/2}+\sin{\tilde{\theta}/2}}\right]}
\end{equation}

\bibliographystyle{myplainnat}
\bibliography{bookbib}

\end{document}